\documentstyle[amsfonts,amssymb]{article}

\newcommand{\msk}{\medskip}
\newcommand{\ndt}{\noindent}

\newcommand{\beq}{\begin{eqnarray}}
\newcommand{\eeq}{\end{eqnarray}}
\newcommand{\beqst}{\begin{eqnarray*}}
\newcommand{\eeqst}{\end{eqnarray*}}

\newcommand{\sgn}{\mbox{ \rm {\small sign}\,}}

\newcommand{\R}{{\mathbb R}}

\newcommand{\qed}{\hfill $\square$}

\newcommand{\dsp}{\displaystyle}

\newtheorem{theorem}{Theorem}[section]
\newtheorem{lemma}[theorem]{Lemma} 
 
\newtheorem{proposition}[theorem]{Proposition}

\newtheorem{remark}[theorem]{Remark}

\title{\bf 
Semilinear Klein-Gordon equation in space-time of   black hole, which is gaining mass in the universe with accelerating expansion}

\author{{\bf Karen Yagdjian and Anahit Galstian   } }

\begin{document}

\date{}

\maketitle

\thispagestyle{empty}
\begin{center}
{\small {School of Mathematical and Statistical Sciences,\\
University of Texas RGV, 
1201 W.~University Drive,  \\
Edinburg, TX 78539,
USA \\
e-mails: karen.yagdjian@utrgv.edu, \hspace{0.2cm} anahit.galstyan@utrgv.edu}}
\end{center}
\medskip
\vspace*{-0.6cm}
  
\addtocounter{section}{0}
\renewcommand{\theequation}{\thesection.\arabic{equation}}
\setcounter{equation}{0}
\pagenumbering{arabic}
\setcounter{page}{1}
\thispagestyle{empty}
\begin{abstract}

In this paper, we consider the propagation of waves in the space-time of a single black hole with a static Schwarzschild 
radius in the expanding universe, namely, the solutions of the linear and semilinear Klein-Gordon equations.

\end{abstract}

\section{Main results}

\setcounter{equation}{0}
\renewcommand{\theequation}{\thesection.\arabic{equation}}

In this paper, we consider the propagation of waves in the space-time of a single black hole (BH) with a static Schwarzschild radius in the expanding universe, namely, the solutions of the linear and semilinear Klein-Gordon equations. 
The line element of this  space-time, which we consider,  is given by 
\begin{equation}
\label{1.2}
\dsp ds^2= - \left( 1-\frac{R_{Sch}}{r} \right) \, c^2dt^2
+  \left( 1-\frac{R_{Sch}}{r} \right)^{-1} a ^2 (t)dr^2 + a^2(t)r^2(d\theta ^2 + \sin^2 \theta \, d\phi ^2),
\end{equation}
where the Schwarzschild radius is 
\[
R_{Sch}:=\frac{2G M_{bh} }{c^2 }, \quad    M_{bh}=constant\,.
\]
Here, $G$ is Newton's gravitational
constant, and $c$ is the vacuum speed of light.
In this paper, we consider  the case of   space-time  with the scale factor $a (t)= \alpha_{acc}  t^{\ell}$, $\ell>1$. 
In what follows, we set $ a_{acc}=1$ and $c=1$, if  no obscurity arises.

Thus, we assume that  the universe is expanding with   acceleration and    that   the mass  of the BHs is changing in time. (See, e.g., \cite{Farrah}, \cite{Perlick}.) 
In \cite{Farrah},   it is suggested to assume that the gaining mass of BH is 
proportional to the scale factor raised to some power $k$. In this paper, we set $k=1$, and consequently, the Schwarzschild radius is static.

In \cite{Macao,ArXiv2023}, the BH with the static   Schwarzschild radius embedded in the 
FLRW space-time with de~Sitter scale factor  $  a(t)=e^{Ht}$  was studied.

We consider another family of the Klein-Gordon equations for the self-interacting waves, with the mass term ingoing redshift proportional to the Hubble parameter $ \dot{a} (t)/ a (t) $, that is, $m(t)=t^{-1} m_c $, where  $m_c:= {m c^2}/{h }$.  For a 
justification of this assumption, see, e.g., \cite{NatarioAHP2022,JDE2021}. Thus, the equation is,
\begin{equation}
\label{1.4} 
    \frac{\partial^2 \psi }{\partial t^2}
+ \frac{3 \ell}{t}    \frac{\partial \psi }{\partial t}
-  t^{  -2\ell}{\mathcal{A}}(x,\partial_x)\psi +   t^{-2}\frac{m^2c^4}{h^2}\left( 1-\frac{R_{Sch}}{r} \right)\psi  
=c^2\left( 1-\frac{R_{Sch}}{r} \right)\Psi(x,\psi)\,, \nonumber
\end{equation} 
where ${\mathcal{A}}(x,\partial_x) $, written in spherical coordinates, is 
\begin{equation}
\label{OpA}
{\mathcal{A}}(x,\partial_x)  
   :=    
 \frac{c^2}{\alpha_{acc}^2}\Bigg\{
    \left( 1-\frac{R_{Sch}}{r} \right)^2 \frac{\partial^2   }{\partial r^2}  
 +\frac{2}{  r   }    \left( 1-\frac{ R_{Sch}}{2r} \right) \left( 1-\frac{R_{Sch}}{r} \right)\frac{\partial   }{\partial r}
 + \left( 1-\frac{R_{Sch}}{r} \right) \frac{1}{r^2}\Delta _{{\mathbb S}^2}   \Bigg\}     \,,  
\end{equation}
and $ \Delta _{{\mathbb S}^2}$ is the Laplace operator on the unit sphere $ {\mathbb S}^2 $. 
According to  Lemma~2.1\cite{ArXiv2023},  the  space-time of the  BH   with a static Schwarzschild radius 
 is the only one  that dissipates   the waves   independently of  their  location in   space. This  is related to the cosmological principle.

We can generalize (\ref{1.4})  as follows
\begin{equation}
\label{KGE_KF}
    \frac{\partial^2 \psi }{\partial t^2}
+ \frac{3 \ell}{t}    \frac{\partial \psi }{\partial t}
-  t^{  -2 \ell}{\mathcal{A}}(x,\partial_x)\psi +    t^{-2}\frac{m^2c^4}{h^2}\psi+V(r,t)  \psi  =c^2\left( 1-\frac{R_{Sch}}{r} \right)\Psi(x,\psi)\,.
\end{equation}
In (\ref{KGE_KF})   $V(r,t) $ is the potential that, in particular,  includes the case of the  gravitational potential 
$
V(r,t)=- t^{-2}\frac{m^2c^4}{h^2} \frac{R_{Sch}}{r}$. 
The term  
$ \Psi(x,\psi)$ represents the self-interaction of the field. 
Consider  the  radial null geodesics  of the space-time (\ref{1.2}) with $\ell >1$, when 
\[
A(t) := \frac{c }{\alpha_{acc}}\left( \frac{t^{1- \ell}}{1- \ell}-\frac{1}{1- \ell} \right) \leq A(\infty)= \frac{c }{\alpha_{acc}}\frac{1}{ \ell-1} \quad \mbox{\rm for all} \quad t \in [1,\infty). 
\]
We consider waves (solutions of the equation) in the exterior  $B^{ext}_{Sch}:=\{x \in \R^3\,|\, |x| > R_{Sch} \}$ of the BH.  The choice of the initial data with supports   sufficiently distant from the horizon of the BH allows us to avoid the manifold,  where the equation degenerates. Hence, we set  
\begin{equation}
\label{13.13}
\mbox{\rm supp} \,\psi _0, \mbox{\rm supp} \,\psi _1 \subseteq \left\{ x \in \R^3\,\left|\, |x|\geq R_{ID}> \frac{c }{\alpha_{acc}}\frac{1}{\ell-1}+ R_{Sch} \right\} \subset  B^{ext}_{Sch}.  \right.
\end{equation}

 Let $H_{(s )}:= H_{(s )}( \R^3)$ be the Sobolev space. The set of all smooth functions with uniformly bounded  derivatives of any order will be denoted ${\mathcal B}^\infty(B^{ext}_{Sch}\times[0,\infty)) $. For the projection operator on $\R^3$ we write $\pi_x$. The space of solutions is
\[
X({R,H_{(s)},\gamma})  := \left\{ \psi  \in C([0,\infty) ; H_{(s)} ) \; \Big|  \;
 \parallel  \psi   \parallel _X := \sup_{t \in [1,\infty) } t^{\gamma}  \parallel
\psi  (x ,t) \parallel _{H_{(s)}}
\le R \right\}\,,
\]
where $\gamma \geq 0$, with the metric
\[
d(\psi _1,\psi _2) := \sup_{t \in [1,\infty) }  t^{\gamma }  \parallel  \psi _1 (x , t) - \psi _2 (x ,t) \parallel _{H_{(s)}}\,.
\]
\smallskip

We  apply   the integral transform approach developed in \cite{Macao}.  (One can refer to \cite{NA2015,Nakamura2014,Costa} on the approach to the semilinear Klein-Gordon equation in the FLRW  spacetimes based on the energy estimate.) Those
integral transforms make it possible to bring the problem with unbounded time down to the problem with
bounded time. This is possible since the space-time has a permanently restricted domain of influence.
\smallskip

Thus, the setting of the   Cauchy Problem is as follows. Let the number $R_{ID} $ be such that $R_{ID} > R_{Sch}$. For every given $\psi_0     \in H_{(s)}   $,  $   \psi_1 \in H_{(s-1)}   $ with supp\,$\psi_0\, \cup$ supp\,$\psi_1  \subseteq \{ x\in  \R^3\,|\, |x| \geq  R_{ID}\} \subset B^{ext}_{Sch}$,  find  
 {\sl a global  in time solution}  $\psi \in C^1([1,\infty); H_{(s)})\cap C ([1,\infty); H_{(s-1)})$  of the equation (\ref{KGE_KF}), such that supp~$\psi (t) \subseteq B^{ext}_{Sch}$ for all times $t >1$ 
and which takes the  initial values
$\psi (x,1) =  \psi_0 (x)$, $ \psi_t (x,1) =  \psi_1 (x)$,   for all $ x \in \R^3$. 
\medskip

 We relate the properties of the 
potential $V(x,t)\in C^2 (B^{ext}_{Sch}\times[0, \infty))$ to the setting 
of the Cauchy problem, more precisely, with   
the support of the initial functions. 
We assume that the potential    $V(x,t)\in {\mathcal B}^\infty (B^{ext}_{Sch}\times[1, \infty))$
 is such that  for given $s$ there are $\varepsilon_P>0 $ and $\delta \geq 0 $ such that  
\begin{equation}
\label{V}
 \|V (x,t) \Phi (x) \| _{H_{(s)}}\leq t^{-\delta}\varepsilon_P \| \Phi   \| _{H_{(s)}}   
 \quad \mbox{\rm for all}\quad  t \in [1,\infty) \quad \mbox{\rm and all} \quad \Phi \in H_{(s)}  \, , 
\end{equation}
with  compact supp~$\Phi  \subseteq B_P \subset B^{ext}_{Sch}$. 
  The number $\varepsilon_P$ may depend on the distance between the sphere $\partial B_P$ and $B_{Sch}$.  
\medskip

 The non-linear term  is supposed to satisfy the following Lipschitz continuity  condition ($\mathcal L$).  
   The smooth in $x \in B^{ext}_{Sch}$ function $\Psi =\Psi (x,\psi )$ is
said to be Lipschitz continuous with exponent $\alpha \geq 0 $ in the space $H_{(s)} $   if supp~$\Psi (x,\psi )\subseteq $ supp~$ \psi $ and there is a constant 
$C_s \geq 0$  such that 
\[
 \|  \Psi (x,\psi  _1 (x))- \Psi (x,\psi  _2(x) ) \|_{H_{(s)} }  \leq C_s
\| \psi  _1 -  \psi _2   \|_{H_{(s)} }
\left( \|  \psi _1  \|^{\alpha} _{H_{(s)} }
+ \|  \psi _2   \|^{\alpha} _{H_{(s)} } \right) \,\, \mbox{ \rm for all} \,\,  \psi  _1, \psi  _2 \in  H_{(s)}\,.
\]
The semilinear terms $\Psi(\psi)=|\psi|^{1+\alpha}$ and $ \Psi(\psi)=\psi|\psi|^{\alpha}$ are of special interest.

\begin{theorem}
\label{T_large_mass}
Consider the Cauchy problem for the equation (\ref{KGE_KF})
in $  \R^3\times [1,\infty) $  with the initial conditions
\begin{equation}
\label{13.12}
\psi (x,1)=\psi _0(x)\in H_{(s)},\quad \partial_t \psi (x,1)=\psi _1(x)\in H_{(s)}    \,,
\end{equation}
with (\ref{13.13}). 
Assume that   the potential $V(x,t)\in {\mathcal B}^\infty (B^{ext}_{Sch}\times[1, \infty))$ satisfies   (\ref{V}) with $\delta=2$\,.
Assume also that the physical mass  of the field is large, that is, 
$
 4m_c^2\geq(3\ell-1)^2$, 
and the nonlinear term $ \Psi(x,\psi)$ satisfies condition $({\mathcal L})$ with $\alpha> \frac{4}{3\ell-1}$ and $\Psi(x,0) =0$. 

If $\varepsilon_P$ and the norms 
 $
 \|\psi _0\|_{H_{(s)}} , \|\psi _1\|_{H_{(s)}}  $ 
with $s> 3/2 $ are sufficiently small, then the Cauchy problem (\ref{KGE_KF})\&(\ref{13.13})\&(\ref{13.12})  has a global solution
$
\psi \in C^2([1,\infty);H_{(s)})$.
 
Moreover, if we denote $M:=\frac{ \sqrt{4 m_c^2-(1-3 \ell )^2}}{ 2( \ell - 1 ) }  $, then  the solution $ \psi$ belongs to the space $X({R,H_{(s)},\gamma})  $, with $ \gamma \leq \ell $  such that: 

\ndt
(a) if $ \alpha >\frac{\ell +3}{2\ell} $, then  $   \gamma> \frac{3 (\ell +1)}{2(\alpha +1)} $ for $M>0$, while $ \ell> \gamma> \frac{3 (\ell +1)}{2(\alpha +1)} $ for $M=0$; \\ 
(b)  if  $\alpha >\frac{4}{3\ell-1} $, then  $ \frac{2}{\alpha}\leq \gamma    \leq \frac{3 (\ell +1)}{2(\alpha +1)}  $ for $M>0$, while  { $ \frac{2}{\alpha}\leq \gamma   \leq \frac{3 (\ell +1)}{2(\alpha +1)}  $ and $\gamma <\ell$}    
for    $M=0$.

Thus, the solution  $ \psi$ decays according to 
$
  \parallel
\psi  (x ,t) \parallel _{H_{(s)}}
\le R t^{-\gamma} , \quad  t \in [1,\infty)$. 
\end{theorem} 

Note, $\dsp \lim_{\ell \to \infty}\frac{\ell +3}{2\ell}=\frac{1}{2 }$ and   $2>  \frac{\ell +3}{2\ell}> \frac{4}{3\ell-1} $ for all $\ell >1$, while $ \frac{\ell +3}{2\ell}= \frac{4}{3\ell-1} $  if $\ell =1$.   
\begin{theorem}
\label{T_small_mass}
Consider the Cauchy problem for the equation (\ref{KGE_KF})
in $  \R^3\times [1,\infty) $  with the initial conditions (\ref{13.12}) 
satisfying (\ref{13.13}). 
Assume that   the potential $V(x,t)\in {\mathcal B}^\infty (B^{ext}_{Sch}\times[1, \infty))$ satisfies   (\ref{V}) with $\delta \geq 2 $\,.
Assume also that the physical mass $m $ of the field is small, that is, 
$4m_c^2 < (3\ell-1)^2 $,  
and the nonlinear term $ \Psi(x,\psi)$ satisfies condition $({\mathcal L})$ with   
$\alpha>\frac{4}{3\ell-1-\sqrt{(1-3 \ell )^2-4 m_c^2}}  $ and $\Psi(x,0) =0$. 

If $\varepsilon_P$ and the norms 
 $
 \|\psi _0\|_{H_{(s)}} , \|\psi _1\|_{H_{(s)}}  $ 
with $s> 3/2 $ are sufficiently small, then the Cauchy problem (\ref{KGE_KF})\&(\ref{13.13})\&(\ref{13.12})  has a global solution
$
\psi \in C^2([1,\infty);H_{(s)})$. 

  Moreover,  if $\delta = 2 $, then the solution $ \psi$ belongs to the space $X({R,H_{(s)},\gamma})  $, with $ \gamma $  such that $\gamma \leq \ell$ and  
\begin{eqnarray*}
&  &\gamma < \frac{3\ell-1-\sqrt{(1-3 \ell )^2-4 m_c^2}}{2 } \,,
\end{eqnarray*} 
 that is, the solution  $ \psi$ decays according to 
$  \parallel
\psi  (x ,t) \parallel _{H_{(s)}}
\le R t^{-\gamma } $, $t \in [1,\infty)$. 

If $\delta > 2 $, then one can choose $\gamma \leq \min\left\{ \ell, \frac{1}{2}\left(3\ell-1-\sqrt{(1-3 \ell )^2-4 m_c^2}\right)  \right\} $.
\end{theorem}

Note, if $\delta > 2 $ and $(3\ell - 1)^2/4 > m^2_c \geq \ell(2\ell-1)$ in the last theorem, then $\gamma$ attains its largest value $\gamma=\ell$.
 
\medskip

The Sobolev embedding  theorem implies the decay of $L_\infty$-norm  of solution. In particular,  we recover some results of \cite{NatarioAHP2022} if the charge of the model is zero. Further, the results 
of  Theorems~\ref{T_large_mass},\ref{T_small_mass}, if $M_{bh}=0 $, agree 
with those  previously known, see  \cite{NA2015,NatarioAHP2023,JMP2019} and references therein. Note, $\lim_{m \to 0} \frac{4}{3\ell-1-\sqrt{(1-3 \ell )^2-4 m_c^2}} =\infty$.

 In \cite{Tsutaya-Wakasugi}, Tsutaya and Wakasugi proved that for solutions of semilinear wave equations, that is,  in the case of ${\mathcal{A}}(x,\partial_x)=\Delta $, $\ell > 1 $, $M_{bh}=0 $, $\Psi(\psi)=|\psi|^{1+\alpha}$,  and $m=0 $, a blow-up in a finite time can
occur for all $\alpha > 0$.

   For the  semilinear Klein-Gordon equation (\ref{1.4}) with ${\mathcal{A}}(x,\partial_x)=\Delta $,  $R_{Sch}=0 $, $V(r,t)=0$,  $\Psi(x,\psi)=|\psi|^{1+\alpha}$, $\alpha>0$,  and integrable speed of propagation $\ell > 1 $, that is,  for the equation with a permanently restricted domain of influence,  to prove a global existence  Aslan, Ebert, and Reissig    \cite{Aslan} applied a partial Fourier transform 
   and a splitting of the phase space into three zones. In particular, Theorem~2.3~\cite{Aslan}  shows  that the condition $\alpha>\frac{4}{3\ell-1-\sqrt{(1-3 \ell )^2-4 m_c^2}}  $ is necessary for the global existence. At the same time the decay rate given by Theorem~\ref{T_small_mass} coincides with one of  Theorem~2.1~\cite{Aslan}.

 Natario and Sasane \cite{NatarioAHP2022}   obtained   exact decay rates for solutions to the Klein-Gordon equation in the de~Sitter universe in flat FLRW form.   For the derivative of the solutions to the massless equation in the cosmological region of the Reissner--Nordstr\"om--de Sitter model, under some assumption on the boundedness of the smooth solution  in the norms over the two components of the future cosmological horizon, parametrized by the flow of the global Killing vector field $\partial/\partial t$, they proved $L_\infty $-norm decay at infinity. In the last model, the mass and charge of the black hole are assumed to be static.

\begin{remark}
 Theorems~\ref{T_large_mass},\ref{T_small_mass} can be extended to the equation 
\[
    \frac{\partial^2 \psi }{\partial t^2}
+ \frac{\mu}{t}    \frac{\partial \psi }{\partial t}
-  t^{  -2 \ell}{\mathcal{A}}(x,\partial_x)\psi +    t^{-2}\frac{m^2c^4}{h^2}\psi+V(r,t)  \psi  =c^2\left( 1-\frac{R_{Sch}}{r} \right)\Psi(x,\psi)\,,
\]
with the parameter $ \mu \in \R$, which is not necessarily relevant to the dimension  of the spatial variable (see, e.g., \cite{Gibbons});   they can be extended to $n$-dimensional black holes. Moreover, to derive a better decay at infinite time, one can appeal to the estimate in homogeneous Sobolev spaces and to the Gagliardo-Nirenberg inequality.  
\end{remark}
\medskip

 The rest of this paper is organized as follows. In Section \ref{S2}, the Liouville transformation and representation of 
the solution to the linear equation are given. Moreover, the equation is transformed into the self-adjoint 
Cauchy-Kowalewski form. The case  of large mass is considered in  Section~\ref{S3}, while the case  of small mass is 
discussed in   Section~\ref{S4}. The proof of existence  of the solution of semilinear equation is based on Banach's fixed point theorem.  The decay estimates for the solution of the linear equation are derived in subsections~\ref{SS3.2}-\ref{SS3.5} and subsections~\ref{SS4.3}-\ref{SS4.6}, respectively. A significant part of the text (subsections~\ref{SS3.1},\ref{SS3.2},\ref{SS4.2},\ref{SS4.3}) is devoted to the estimation of the integrals of the kernel functions of the integral transform containing   Gauss's hypergeometric function.  

\section{Preliminaries. Linear equation}

\label{S2}

\setcounter{equation}{0}
\renewcommand{\theequation}{\thesection.\arabic{equation}}

We will use the notations $\vec{ x}=(x_1,x_2,x_3):=(x,y,z)$, $x_0=t$, and $\vec{\xi }=(\xi _1,\xi _2,\xi _3)$. The inner   product in $\R^3$ will be denoted $\vec{ x}\cdot \vec{\xi }$. 
In  the case of $\ell=0$, the linear Klein-Gordon equation without source in Cartesian coordinates can be written (see, e.g., \cite{ArXiv2023}) in the form  (\ref{1.4}), 
where
\begin{eqnarray}
\label{SymbA}
{\mathcal A}(\vec{ x};\vec{\xi } )
 & = &
  A_2(\vec{ x};\vec{\xi } )+A_1(\vec{ x};\vec{\xi } )\,,
\end{eqnarray}
while $ A_2(\vec{ x};\vec{\xi } )$ and $A_1(\vec{ x};\vec{\xi } )$ are the principal symbol and the low order symbol, respectively, and
\[
F(|\vec{x}|)=F(r):=  1-\frac{R_{Sch}}{ |\vec{x}|},\quad r:=|\vec{x}|:=\sqrt{x^2+y^2+z^2}, \quad r >R_{Sch}\,.
\]
Thus, the symbol $ {\mathcal A}(\vec{ x};\vec{\xi } )$ of the operator $ {\mathcal A}(\vec{ x},\partial_x)$ (\ref{OpA})  is given by  
\[
 A_2(\vec{ x};\vec{\xi } ) 
   =    
-c^{2}F(r)\Bigg( | \vec{\xi }|^2 -\frac{ R_{Sch}  \left(\vec{ x}\cdot \vec{\xi }\right)^2}{|\vec{ x}|^3 }\Bigg) ,\quad
A_1(\vec{ x};\vec{\xi } ) 
   =   
-c^{2}F(r)\frac{  i R_{Sch}  \left(\vec{ x}\cdot \vec{\xi }\right)}{|\vec{ x}|^3 }.
\]

On every compact set in $B_{Sch}^{ext}\times \R $, the equation (\ref{KGE_KF}) is strictly hyperbolic. It 
 has multiple characteristics on the event horizon $r=R_{Sch}$. 
Consequently, the well-posedness of the Cauchy problem requires some kind of Levi condition. (See, for details,
\cite{YagBook}.) \\

\noindent{\bf The finite propagation speed property. } 
Consider the equation
 \begin{equation}
 \label{2.2}
    \frac{\partial^2 \psi }{\partial t^2}
+ \frac{3 \ell}{t}    \frac{\partial \psi }{\partial t}
-  t^{-2\ell}{\mathcal{A}}(x,\partial_x)\psi +  t^{ - 2 } F (r) m_c^2  \psi = 0\,,
\end{equation}
in  $B^{ext}_{Sch}$. 
We apply the  Liouville transform with $k_\pm$ defined in subsection~\ref{SS2.2} 
$
\psi  (x,t)= t^{k_\pm} \sqrt{F(r)}u (x,t)
$
to (\ref{2.2}), then the covariant   Klein-Gordon  equation (\ref{2.2}) reads 
\begin{eqnarray} 
\label{2.2b}
 &  &
  \frac{1}{c^2}\frac{\partial^2   u }{\partial t^2} 
-  t^{-2\ell}{\mathcal A}_{3/2}(x,\partial_x)u 
+    t^{-1}(\ell+2 i {\cal{M}})  \frac{\partial   }{\partial t } u   -\frac{R_{Sch}}{r}    t^{-2}\frac{m_c^2  }{c^2 }  u   =0\,, 
\end{eqnarray}
where $ {\cal{M}} $ is written in (\ref{Mpm}) and  the  operator ${\mathcal A}_{3/2}$ in the spherical coordinates is defined (see, e.g.,  \cite[Sec.2.3]{ArXiv2023}) by
\[ 
{\mathcal A}_{3/2}(x,\partial_x)v 
   :=  
    F (r)^{3/2} \frac{\partial^2    }{\partial r^2}\sqrt{F(r)}v
+\sqrt{F(r)}\frac{2}{  r   }    \left( 1-\frac{R_{Sch}}{2r} \right)  \frac{\partial   }{\partial r}   \sqrt{F(r)}v
+ F(r) \frac{1}{  r^2 } \Delta_{S^2} v \,,\nonumber 
\]
while the term
$
-t^{-2}\frac{R_{Sch}}{r}   \frac{m_c^2  }{c^2 }
$
 can be regarded as a potential $V=V( \vec{x} ,t)$. 
According to  \cite{ArXiv2023},
the operator ${\mathcal A}_{3/2}(x,\partial_x)$   can be written as follows
\[
{\mathcal A}_{3/2}(\vec{x};D_x)v=\sum_{i,j=1}^3\frac{\partial}{\partial x_i}\left( a_{ij}(\vec{x})\frac{\partial}{\partial x_j}v \right)  +\frac{G^2 M^2_{bh}}{c^4  |\vec{x}|^4}v\,,
\]
with the coefficients $ a_{ij}(\vec{x}) $   such that
\begin{eqnarray*}
a_{ii}(\vec{x}) 
& := &
-F(r) \left(\frac{R_{Sch} x_i^2}{ |\vec{x}|^{3}}-2\right), \quad i=1,2,3,\\
a_{ji}(\vec{x}) 
& = & 
a_{ij}(\vec{x}):= 
-F(r)\frac{ R_{Sch} x_i x_j  }{ |\vec{x}|^{3}},\quad i,j=1,2,3, \quad i\not=j\,.
\end{eqnarray*}
The finite propagation speed property for the linear equation (\ref{2.2b}) can be proved by standard arguments (see, e.g.,\cite{Taylor}).
 \medskip
 
\noindent
{\bf The radial geodesics.}  
Next, consider   the radial geodesic solving 
\begin{equation}
\label{r_geodesic}
\frac{dr}{dt}=-  \frac{c}{\alpha_{acc} }  t^{-\ell}  \left( 1-\frac{R_{Sch}}{r} \right)  
\end{equation}
and starting at $R_{ID}$. The existence of global in time  geodesic is given by the following statement.
\begin{lemma}
The implicit function   $r=r(t)$, which is  defined for all $t \in [1,\infty)$ by 
\begin{equation}
\label{r_geod_form}
 R_{ID}-r -R_{Sch}  \ln \left(   \frac{  r-R_{Sch}}{R_{ID}-R_{Sch}}  \right)   =\frac{c }{\alpha_{acc}}\left( \frac{t^{1- \ell}}{1- \ell}-\frac{1}{1- \ell} \right),
\end{equation}  
with some positive number $\varepsilon$ satisfies the inequality
\begin{equation}
\label{r_support}
r(t) \geq R_{Sch}+\varepsilon \quad \mbox{\it for all} \quad t >1\,.
\end{equation}
\end{lemma}

\ndt
{\bf Proof.}  The solution to (\ref{r_geodesic}) is given by (\ref{r_geod_form}). 
For $\tau \in [0,  \infty)  $, consider an implicit function $z(\tau)$ given by the problem
\[
z -R_{Sch}  \ln \left( 1- \frac{ z}{R_{ID}-R_{Sch}}  \right)   =\tau,\quad  z(1)=0.
\]   
The function $z=z(\tau)$ is well defined if $z(\tau) \in [0, R_{ID}-R_{Sch})$. Moreover, it is positive, and is continuous. Indeed, 
\[
\frac{dz}{d\tau}\left(1 + R_{Sch}   \frac{ 1}{R_{ID}-R_{Sch}-z }   \right)   =1
\]
as long as $z(\tau)< R_{ID}-R_{Sch}$ implies $\frac{dz}{d\tau}>0 $. Denote $z_\infty=z( A(\infty))$, where $A(\infty)= \frac{c}{\alpha ({\ell}-1)}$, then 
$z(\tau) \leq z_\infty$ for all  $\tau \in [0,cA(\infty)]$.   
The existence of the number $z_\infty\in [0, R_{ID}-R_{Sch})$ follows from 
\[
\lim_{z  \nearrow  R_{ID}-R_{Sch} } \left(  z -R_{Sch}  \ln \left( 1- \frac{ z}{R_{ID}-R_{Sch}}  \right) \right)  = \infty \,.
\] 
Hence, there is a number $\varepsilon>0$ such that 
$
z(\tau) \leq R_{ID}-R_{Sch} -\varepsilon$  for all $\tau \in [0,cA(\infty)]$.  
Next, we define $\tau$   by
\[
 \tau (t):=A(t)\quad \mbox{\rm and} \quad \frac{d\tau}{dt}= \frac{c}{a(t)} >0 \quad \mbox{\rm for all} \quad t \in [1,\infty)\,.
\]
Thus, the function $r=r(t)$, which is  defined on $[0,\infty) $ and  
$r(t)= R_{ID}-z(\tau (t))$,   where $\tau (t)< cA(\infty)$  for all $t \in [1,\infty)$. 
The lemma is proved. \qed 
\smallskip

Hence, if the initial functions have compact supports  in $K\subset  B_{Sch}^{ext}\subset \R^3$, dist$(\partial K,B_{Sch})>\frac{c}{\alpha _{acc} (\ell -1)}$, then there is $\varepsilon>0$ such that for all  $t \in [1,\infty)$, the solution  of (\ref{KGE_KF}) with $\Psi=0$  has a compact support in $ \{x \in \R^3 \,|\, |x| > R_{Sch}+\varepsilon\} \subset B_{Sch}^{ext}$.   
This follows from the dependence domain. 
 \smallskip 

Thus, for every given fixed compact $K\subset  B_{Sch}^{ext}$, the coefficients of the equation can be continued smoothly outside of the small $\delta$-neighborhood of $K$ to be constants. This continuation does change   solutions, which takes   initial data supported by $K$, at least for time duration $d/c$. Consequently, we can replace the operator ${\mathcal{A}}(x,\partial_x) $ with such continuation (auxiliary operator) 
${\mathcal A}_{\varepsilon }(\vec{ x};\vec{\xi } )$ that has a symbol 
\begin{equation}
\label{AUX} 
{\mathcal A}_{\varepsilon }(\vec{ x};\vec{\xi } )  
  :=  
c^{2}\left(1-\chi_\varepsilon (\vec{ x})\frac{R_{Sch}}{|\vec{ x}|}\right)\left(-|\xi |^2    +\chi_\varepsilon (\vec{ x})\frac{ R_{Sch} \left(\vec{ x}\cdot \vec{\xi }\right)^2}{  |\vec{ x}|^{3} } 
-\chi_\varepsilon (\vec{ x})\frac{ i R_{Sch}  \left(\vec{ x}\cdot \vec{\xi }\right)}{   |\vec{ x}|^{3} }\right),  
\end{equation}
where  $\chi_\varepsilon $ is a cutoff function such that $\chi_\varepsilon (\vec{ x})=0$ if $|\vec{ x}| < R_{Sch}+\varepsilon/2 $ while $\chi_\varepsilon (\vec{ x})=1$ if $|\vec{ x}| > R_{Sch}+\varepsilon  $. 
This allows us to apply the $H_{(s)}$-estimates.

\medskip

\subsection{The partial Liouville transform}
\label{SS2.2}

To get rid of the mass term in the initial value  problem  for the equation 
\begin{eqnarray*}
\frac{\partial^2 }{\partial t^2} \psi+\frac{3\ell}{t}  \frac{\partial }{\partial t}\psi - t^{-2\ell}{\mathcal{A}}(x,\partial_x)\psi +t^{ - 2 }m_c^2 \psi =f 
\end{eqnarray*}
we use the partial Liouville transform $\psi  (x,t)= \omega (t) u (x,t)$.  
Then the equation reads
\[
u_{tt} (x,t)- t^{-2\ell}{\mathcal{A}}(x,\partial_x)  u(x,t) 
+\left( 2\frac{\omega_t (t)}{\omega (t)}  +\frac{3\ell}{t}  \right)  u_t (x,t)+\left(\frac{\omega_{tt}(t)}{\omega (t)} + \frac{3\ell}{t}  \frac{\omega_t (t)}{\omega (t)} +\frac{ m_c^2}{t^{2 }}  \right) u (x,t)  =\frac{1}{\omega (t)}f\,. 
\]
We set
\begin{eqnarray*}
 \frac{\omega_{tt}(t)}{\omega (t)} + \frac{3\ell}{t}  \frac{\omega_t (t)}{\omega (t)}+ \frac{ m_c^2}{t^{2 }} =0\,.
\end{eqnarray*}
If we look for
$
\omega (t)=t^k,
$
then 
$k(k-1)+3\ell k+m_c^2=0$. 
The solutions to the last   equation are
\begin{equation}
\label{kpm}
k_{\pm}=\frac{1}{2} \left(1-3 \ell  \pm\sqrt{(1-3\ell  )^2 -4 m_c^2}\right).
\end{equation}  
We consider the equation with $\ell>1$ that implies that we have to consider three cases: 
(1)  $m_c<\frac{3\ell-1}{2}$, 
(2)  $m_c=\frac{3\ell-1}{2}$  $\Longrightarrow  k_{\pm}=\frac{1}{2} (1-3 \ell )$, and 
(3)  $m_c>\frac{3\ell-1}{2}$. 
We say that {\it mass is small} if $m_c<\frac{3\ell-1}{2}$, while in other cases we say mass is  {\it large}. 
\bigskip

Thus,  we arrived at  the generalized Euler-Poisson-Darboux  equation  
\begin{eqnarray}
\label{WEDT}
u_{tt} (x,t)- t^{-2\ell}{\mathcal{A}}(x,\partial_x)  u(x,t) +  \frac{1\pm\sqrt{(1-3 \ell )^2-4 m_c^2}}{t}   u_t (x,t) =t^{-k_\pm}f \,.
\end{eqnarray}
Denote 
\begin{eqnarray}
\label{Mpm}
&  &
  {\cal{M}}_\pm :=\frac{1}{2 i}( 1-\ell \pm\sqrt{(1-3 \ell )^2-4 m_c^2} ),  
\quad \frac{1\pm\sqrt{(1-3 \ell )^2-4 m_c^2}}{t} := t^{-1}(\ell+2 i {\cal{M}_\pm}) \,.
\end{eqnarray}
Then, for the case of small mass, that is,  $4 m_c^2<(1-3 \ell )^2 $, we have
$ 2i \frac{{\cal{M}_-}}{ 1-\ell }  >1$, $ 2i \frac{{\cal{M}_+}}{ 1-\ell }  <1$, 
while for the large mass, that is,  $4m_c^2 \geq(1-3 \ell )^2   $, we have
$
\Re  2i \frac{{\cal{M}_-}}{ 1-\ell } =\Re  2i \frac{{\cal{M}_+}}{ 1-\ell } =1\,.
$

\subsection{Representation of the solution of linear equation}

We set the initial time $\varepsilon =1$ in the formulas of \cite{JDE2021} (for some particular cases, see also \cite{NatarioAHP2023,Palmieri}), then
\begin{eqnarray*}
&  &
\phi (1) :=\frac{1}{1-\ell}  <0, \quad  \phi (t) :=\phi (1)t^{1-\ell}  \quad 
\mbox{\rm for all} \quad t\geq 1\,,\quad 
 \phi (t)- \phi (b ) >0, \quad \mbox{\rm for all} \quad  1\leq b<t\,.
\end{eqnarray*}
Next, we need the following kernel functions $E(r,t;t_0  ; {\cal{M}}) $, $K_1 \left(r,t; {\cal{M}} ;1 \right) $,  and $K_0\left(r,t  ; {\cal{M}}; 1 \right) $, from \cite{JDE2021}
\begin{eqnarray}
\label{Edef}
 E(r,t;t_0  ; {\cal{M}}) 
& = & 
{ 2^{ \frac{2 i{\cal{M}}}{ 1-\ell }}} (1-\ell )^{\frac{\ell }{1-\ell }}   \phi (t_0  )^{\frac{ \ell+2i{\cal{M}}}{1-\ell }}
\left( \left(\phi (t)+ \phi (t_0  )  \right)^2- r^2 \right)^{-  \frac{ i{\cal{M}}}{ 1-\ell }} \\
&  &
\times  F \left(  i \frac{{\cal{M}}}{ 1-\ell },  i\frac{{\cal{M}}}{ 1-\ell };1;\frac{\left( \phi (t)-  \phi (t_0  )   \right)^2-r^2}{\left( \phi (t)+ \phi (t_0  ) \right)^2-r^2}\right)  \,,  \nonumber \\
\label{K1def} 
 K_1 \left(r,t; {\cal{M}} ;1 \right)
& := &
 2^{2 i  \frac{{\cal{M}}}{1-\ell}} \phi ( 1 )^{2 i  \frac{{\cal{M}}}{1-\ell}}   \left( \left(\phi (t)+ \phi (1 )  \right)^2-  r ^2\right)^{-  i \frac{{\cal{M}}}{1-\ell}}    \\ 
&  &
\times F \left(  i \frac{{\cal{M}}}{1-\ell},  i \frac{{\cal{M}}}{1-\ell};1;\frac{\left(\phi (t)- \phi (1 )\right)^2- r ^2}{\left(\phi (t)+ \phi (1 )   \right)^2-  r ^2 }\right)  \,, \nonumber\\ 
\label{K0def}  
K_0\left(r,t  ; {\cal{M}}; 1 \right)
& := &
- i { 2^{1+2 i  \frac{{\cal{M}}}{1-\ell }}}  \frac{{\cal{M}}}{1-\ell }  \phi^{  \frac{2i{\cal{M}}}{1-\ell }} (1 )
\left(\left(\phi (t)+ \phi (1 )\right)^2-r^2\right)^{-  i\frac{{\cal{M}}}{1-\ell }}  \\
&  &
\times \left\{\frac{ r^2- \phi  (t)(\phi  (t)- \phi (1 ))    }{ r^2-\left(\phi (t)- \phi (1)\right)^2 } 
F \left(  i \frac{{\cal{M}}}{1-\ell },  i \frac{{\cal{M}}}{1-\ell };1;\frac{\left( \phi (t)- \phi (1) \right) ^2-r^2}{(  \phi (t)+ \phi (1 ) )^2-r^2}\right)\right.\nonumber\\
&  &
\left.-\frac{2 \phi (t)  \phi (1)   \left(\phi^2 (t) - \phi^2 (1)  -r^2\right)  }{\left(\left(\phi (t)- \phi (1 )\right)^2- r^2\right) \left(\left(  \phi (t)+ \phi (1 ) \right)^2-r^2\right)}\right.\nonumber\\
&  &
\left.\times F \left(  i \frac{{\cal{M}}}{1-\ell }+1,  i \frac{{\cal{M}}}{1-\ell };1;\frac{\left( \phi (t)- \phi (1 ) \right) ^2-r^2}{(  \phi (t)+ \phi (1 ) )^2-r^2}\right)\right\}\,, \nonumber
\end{eqnarray} 
 respectively.\footnote{The definition of $E$ and $K_1$ in   \cite{JDE2021} differs from the one used in this paper by factors.} Here $F(a,b;c;\zeta  ) $ is the  Gauss's
hypergeometric function (see, e.g., \cite{B-E}). We also define  
\begin{eqnarray} 
\label{K2def}
K_2 \left(r,t  ;{\cal{M}};1 \right) := K_0\left(r,t  ;{\cal{M}};1 \right)+  2 i  \frac{\cal{M}}{1-\ell}    K_1 
\left(r,t; {\cal{M}} ;1 \right) \,.
\end{eqnarray} 
We rewrite equation (\ref{WEDT}) as follows:
\begin{eqnarray*} 
  u_{tt} - t^{-2\ell } A(x,\partial_x)u +t^{-1}(\ell+2 i {\cal{M}_\pm})  u_t =  t^{-k_\pm}f (x,t)\,, \quad t>1 \,,\quad x \in \Omega \,,
\end{eqnarray*}
where (\ref{Mpm})   has been used.
\begin{theorem}
\label{T2.4} 
\cite[Theorem 4.1]{JDE2021}
Assume that the function $v=v_g  \left(x,r;b  \right) \in C_{x,r,b}^{d,2,0}$ solves the problem
\begin{eqnarray*}  
\cases{ 
\dsp  v _{r r } (x,r;b)-     A(x,\partial_x) v (x,r;b)   =0\,, \quad x \in \Omega,\cr 
 v(x,0;b )=g(x,b), \quad  v _r (x,0;b )= 0\,,\quad x \in \Omega,} 
\end{eqnarray*}
while  the function $v  =v_\varphi   (x,r )  \in C_{x,r}^{d,2}$ solves the problem
\begin{eqnarray*}  
\cases{ 
\dsp   v _{r r } (x,r )-     A(x,\partial_x) v (x,r )   =0\,,\quad x \in \Omega,\cr 
 v(x,0 )=\varphi (x ), \quad  v  _r   (x,0  )= 0\,,\quad x \in \Omega.} 
\end{eqnarray*}
Then the function $u=u(x,t) $  defined by
\begin{eqnarray*}
  u(x,t )
& = &
 \int_1^{t  } \,d b\int_0^{ \phi (t)- \phi (b ) }  
 E(r,t;b ; {\cal{M}}) v_g  \left(x,r;b  \right) \,dr\\
&  &
+\int_0^{\phi (t)- \phi (1 )}   K_1 \left(r,t; {\cal{M}} ;1 \right) v_{u _1}(x,r)\,dr\\
&   &
+
\left(\frac{\phi (t)}{ \phi (1)} \right)^{-   i\frac{{\cal{M}}}{1-\ell} } 
v_{u _0}\left(x,\phi (t)- \phi (1) \right)
+\frac{1}{\phi (1)}\int_0^{\phi (t)- \phi (1)} K_2 \left(r,t  ;{\cal{M}};1 \right)  
v_{u _0}(x, r)\,dr  \,,
\end{eqnarray*} 
where $\phi (t):=\frac{1}{1-\ell}t^{1-\ell} $ and   the kernels $ E(r,t;b ; {\cal{M}})$, $ K_1 \left(r,t; {\cal{M}} ;1 \right) $, $K_0\left(r,t  ; {\cal{M}}; 1 \right) $, and $K_2 \left(r,t  ;{\cal{M}};1 \right) $  are defined by (\ref{Edef}), (\ref{K1def}), (\ref{K0def}), and (\ref{K2def}), respectively, 
is a solution to the Cauchy problem
\begin{eqnarray*}
\cases{ 
  u_{tt} - t^{-2\ell } A(x,\partial_x)u +t^{-1}(\ell+2 i {\cal{M}})  u_t =  g (x,t)\,, \quad t>1 \,, \quad x \in \Omega\,,\cr
 u(x,1   )=u _0(x), \quad  u _t (x,1   )= u _1 (x),\quad x \in \Omega\,.} 
\end{eqnarray*}
\end{theorem}

We study the Cauchy problem    (\ref{KGE_KF})\&(\ref{13.12})
    through the integral equation.
To define that integral equation we  appeal to the operator 
\begin{equation}
\label{G}
G:={\mathcal K}\circ {\mathcal EE}\,,
\end{equation}
where ${\mathcal EE}$ stands for the evolution (wave) equation in the exterior of BH in the universe without expansion as follows. For the function $f(x,t) $ we define
\[
v(x,t;b):= {\mathcal EE} [f](x,t;b)\,,
\]
where the function 
$v(x,t;b)$   
is a solution to the Cauchy problem 
\begin{eqnarray}
\label{2.15}
&   &
\cases{ \partial_t^2 v - {\mathcal A}(x,\partial_x)v =0, \quad x \in B^{ext}_{Sch} \subset {\mathbb R}^3, \quad t \geq 0, \cr 
v(x,0;b)=f(x,b)\,, \quad v_t(x,0;b)= 0\,, \quad x \in B^{ext}_{Sch} \subset {\mathbb R}^3\,,  \quad b \geq 1,}
\end{eqnarray} 
while ${\mathcal K}$ is introduced  by 
\begin{eqnarray}
\label{9.9}
{\mathcal K}[v]  (x,t) 
&  :=  &
 2t^{k }\int_1 ^{t  } b^{-k } \,d b\int_0^{ \phi (t)- \phi (b ) }  
 E(r,t;b ; {\cal{M}}) v \left(x,r;b  \right) \,dr \,.
\end{eqnarray}
The kernel $ E(r,t; 0,b;M) $ is given by  (\ref{Edef}).

\subsection{Energy estimates and   energy conservation}

Let $\partial_t^2- A(x,\partial_x) =\partial_t^2- \sum_{|\alpha |\leq 2} a_{\alpha }(x)\partial _x^\alpha $ be a  second-order strictly hyperbolic   operator with coefficients $a_{\alpha } \in {\mathcal B}^\infty $, where   ${\mathcal B}^\infty $  is the space of all $C^\infty ({\mathbb R}^3)$ functions with uniformly bounded derivatives of all orders.  
 Let $v=v(x,t) $ be the solution of the problem
\begin{eqnarray*}
&   &
\cases{ \partial_t^2 v - A(x,\partial_x)v =0, \quad x \in {\mathbb R}^n, \quad t \geq 0, \cr 
v(x,0) = v_0 (x),\quad   v_t (x,0)= v_1(x),\quad x \in {\mathbb R}^n\,.}
\end{eqnarray*}
 The following energy estimate is well known. (See, e.g., \cite{Taylor}.) For every $s \in {\mathbb R}$ and given $T>0$ there is  $ C_s(T)$ such that  
\begin{eqnarray}
\label{Blplq}
&   &
 \|v_{t}(t)\|_{H_{(s)}}   +  \|v (t)\|_{H_{(s+1)}} \leq  C_s(T) ( \|v_1\|_{H_{(s)}} +  \|v_0 \|_{H_{(s+1)}} )  , \quad 0 \leq  t \leq T  \,.
\end{eqnarray}
We note that although in the last estimate the time interval is bounded, however, due to the integral transform approach,  it is possible to reduce the problem with infinite time to the problem with finite time and to apply (\ref{Blplq}).
 Here, the permanently restricted domain of influence of the
space-time  with $\ell>1$ comes into play.

We are going to apply the estimate (\ref{Blplq}) to the problem     
\begin{eqnarray*}
&   &
\cases{ \partial_t^2 v - {\mathcal A}_{\varepsilon }(x,\partial_x)v =0, \quad x \in {\mathbb R}^3, \quad t \geq 0,  \cr
v(x,0) = v_0 (x),\quad   v_t (x,0)= v_1(x),\quad x \in {\mathbb R}^3\,,}
\end{eqnarray*}
where the operator ${\mathcal A}_{\varepsilon }(x,\partial_x)$ has a symbol ${\mathcal A}_{\varepsilon }(x;\xi)$ of   (\ref{AUX}). In that case the constant $C_s(T) $ depends on $ \varepsilon $ as well. 
The conservation of the energy of the solution of the equation 
\[
    \frac{\partial^2 \psi }{\partial t^2}
-{\mathcal{A}}(x,\partial_x)\psi +  F(r) \frac{m^2 c^4 }{h^2}  \psi  =0 
\]
is known (see, e.g., \cite{ArXiv2023}).

 \medskip
 
 \subsection{Equation in self-adjoint  form}

The semilinear Klein-Gordon equation without potential is  
\begin{eqnarray}
\label{KGSLKF}
    \frac{\partial^2 \psi }{\partial t^2}
+ \frac{3 \ell}{t}    \frac{\partial \psi }{\partial t}
-  t^{-2\ell}{\mathcal{A}}(x,\partial_x)\psi +  t^{ - 2 } F (r) m_c^2   \psi = c^2F (r)\Psi(\psi)\,,
\end{eqnarray}
where $ {\mathcal{A}}(x,\partial_x)$ is defined in (\ref{OpA}). 
The   operator $ {\mathcal{A}}(x,\partial_x)={\mathcal A}(\vec{ x},\partial_x)={\mathcal{A}}(x,y,z;D_x,D_y,D_z)$ has the symbol (\ref{SymbA}) 
in Cartesian coordinates. 
We apply the  Liouville substitution
$\psi  (x,t)= t^{k_\pm} \sqrt{F(r)}u (x,t)$
to (\ref{KGSLKF}), then the covariant   Klein-Gordon  equation (\ref{KGSLKF}) became 
\begin{eqnarray*} 
  \frac{1}{c^2}\frac{\partial^2   u }{\partial t^2} 
-  t^{-2\ell}{\mathcal A}_{3/2}(x,\partial_x)u 
+   t^{-1}(\ell+2 i {\cal{M}_\pm}) \frac{\partial   }{\partial t } u   -\frac{R_{Sch}}{r}    t^{-2}\frac{m_c^2  }{c^2 }  u   
  =    t^{-k_+}\sqrt{F(r)}\Psi \left(t^{ k_+}\sqrt{F(r)}u \right)\,, 
\end{eqnarray*}
where  the  operator ${\mathcal A}_{3/2}(x,\partial_x)$ in the spherical coordinates is defined (see, \cite[Sec.2.3]{ArXiv2023}) by
\begin{eqnarray*}
{\mathcal A}_{3/2}(x,\partial_x)v 
& := & 
    F (r)^{3/2} \frac{\partial^2    }{\partial r^2}\sqrt{F(r)}v
+\sqrt{F(r)}\frac{2}{  r   }    \left( 1-\frac{ R_{Sch}}{2r} \right)  \frac{\partial   }{\partial r}   \sqrt{F(r)}v
+ F(r) \frac{1}{  r^2 } \Delta_{S^2} v \,,\nonumber 
\end{eqnarray*}
while the term
$
-t^{-2}\frac{R_{Sch}}{r}   \frac{m^2 c^2 }{h^2}
$
 can be regarded as a potential $V=V(x,y,z,t)$.

The symbol ${\mathcal A}_{3/2}(x,\xi ) $  of operator ${\mathcal A}_{3/2} (x,\partial_x) $ in Cartesian coordinates   is
\begin{eqnarray*} 
{\mathcal A}_{3/2}(\vec{x},\vec{\xi}) 
& := & 
\left(1-\frac{R_{Sch}}{ |\vec{x}|}\right) \left(-|\vec{\xi}|^2 +\frac{R_{Sch}\left(\vec{x}\cdot\vec{\xi}\right)^2}{ |\vec{x}|^{3}}\right)
+\frac{R_{Sch}^2}{4 |\vec{x}|^4} .
\end{eqnarray*}
According to Lemma~2.4~\cite{ArXiv2023},
the operator ${\mathcal A}_{3/2}(x,\partial_x)$  is essentially self-adjoint   on $C_0^\infty(B^{ext}_{Sch})$, and on every closed subset in $B^{ext}_{Sch}$ it is an elliptic  operator that is non-positive on the subspace of functions with   supports in $B^{ext}_{Sch}$.

\section{The case of large mass. Proof of Theorem~\ref{T_large_mass}}
\label{S3}

In this case, for the partial Liouville transform we choose $k_-$ and ${\cal{M}_-} $ from (\ref{kpm}) and (\ref{Mpm}), respectively.  
 We also denote
\[
M:=\frac{ \sqrt{4 m_c^2-(1-3 \ell )^2}}{ 2( \ell - 1 ) }\geq 0,\quad \mbox{then} \quad i \frac{{\cal{M}_-}}{ 1-\ell }
=\frac{1}{ 2}+  iM  \,.
\]
In this section, we drop the subindex of $ \cal{M}_-$.
\begin{lemma} 
For $M \in [0,\infty)$  and the given above parameter $\ell $  with some constant $C$ 
\begin{eqnarray*}
 \left| F \left(  \frac{1}{2},   \frac{1}{ 2}+  iM,\frac{3}{2};\frac{\left(z-1\right)^2 }{\left(z+ 1   \right)^2  }\right)\right|  
& \leq & 
C \quad \mbox{  for all}\quad z \in [1,\infty)\,, \\
\left|F \left(\frac{1}{2}, \frac{3}{ 2}+  iM;\frac{3}{2};\frac{(z-1)^2}{(z+1)^2}\right) \right| 
 &  \leq &
C +C(1-\sgn (M))z^{\frac{1}{2}} \quad \mbox{  for all}\quad z \in [1,\infty)\,.   
 \end{eqnarray*}
\end{lemma}
\medskip

\ndt
{\bf Proof.} According to  \cite[(1) Sec.~2.10]{B-E} 
\begin{eqnarray*}  
   F \left(  \frac{1}{2}, \frac{1}{ 2}+  iM,\frac{3}{2};\zeta  \right)  
& = &
 \frac{\Gamma(\frac{3}{2})\Gamma(\frac{1}{ 2}-  iM)}{\Gamma(1)\Gamma(1-iM)}
F \left(  \frac{1}{2}, \frac{1}{ 2}+  iM, \frac{1}{2}+iM;1-\zeta  \right) \\
&  &
+ \left(1 - \zeta  \right)^{\frac{1}{ 2}-iM}\frac{\Gamma(\frac{3}{2})\Gamma(\frac{1}{ 2}-  iM)}{\Gamma(\frac{1}{ 2})\Gamma(\frac{1}{ 2}+iM)}
F \left(  1, 1-  iM, \frac{3}{2}-iM;1-\zeta \right),
\end{eqnarray*} 
where $\zeta =\frac{\left(z-1\right)^2 }{\left(z+ 1   \right)^2  }$ and $1-\zeta = \frac{4 z}{(z+1)^2} \to 0$   as $z \to \infty$. 
Consequently, the first estimate follows. 
Similarly,   with $M> 0$ we obtain 
\begin{eqnarray*} 
\left|F \left(\frac{1}{2}, \frac{3}{ 2}+  iM;\frac{3}{2};\frac{(z-1)^2}{(z+1)^2}\right) \right| 
 &  \leq &
 C   \quad \mbox{\rm for all}\quad z \in [1,\infty),
 \end{eqnarray*}
 while   according to  \cite[(23) Sec.~2.1]{B-E}  with $M= 0$ we obtain 
\begin{eqnarray*} 
\left|F \left(\frac{1}{2}, \frac{3}{ 2}+  iM;\frac{3}{2};\frac{(z-1)^2}{(z+1)^2}\right) \right| 
 &  \leq &
C  + C(1-\sgn (M))\left(   \frac{4z }{ (z+ 1  )^2} \right)^{-\frac{1}{2}} \,.
 \end{eqnarray*}
The lemma is proved. \qed

\begin{lemma}
\label{L3.2}
Let $4m_c^2\geq (3\ell-1)^2   $, then
\begin{eqnarray*}
&    &
  \int_0^{z-1}   \left|    \left( \left(1+ z \right)^2-  y ^2\right)^{-  \frac{1}{ 2} }  F \left(  \frac{1}{ 2}+  iM, \frac{1}{ 2}+  iM;1;\frac{\left(z- 1\right)^2- y ^2}{\left(z+ 1  \right)^2-  y ^2 }\right) \right|  \,dy\\
& \lesssim &
(1+\ln (z))^{1-\sgn M} \quad \mbox{ for all}\quad z \in [1,\infty)\,.
\end{eqnarray*}
\end{lemma}
\medskip

\ndt
{\bf Proof.} 
According to   \cite[(7.3)]{Yag_Galst_CMP} 
\begin{eqnarray}
\label{3.15}
\left| F \left(  \frac{1}{ 2}+  iM, \frac{1}{ 2}+  iM;1;\zeta \right) \right| \leq C(1-\ln (1-\zeta))^{1-\sgn M}
\end{eqnarray}
we obtain 
\begin{eqnarray*}
&  &
\int_0^{z-1}   \left|    \left( \left(1+ z \right)^2-  y ^2\right)^{-  \frac{1}{ 2} }   F \left(  \frac{1}{ 2}+  iM, \frac{1}{ 2}+  iM;1;\frac{\left(z- 1\right)^2- y ^2}{\left(z+ 1  \right)^2-  y ^2 }\right) \right|  \,dy\\
& \leq &
C\int_0^{z-1}   \left|    \left( \left(1+ z \right)^2-  y ^2\right)^{-  \frac{1}{ 2} }  \left(1-\ln \left( \frac{4z}{\left(z+ 1  \right)^2-  y ^2 }\right) \right)^{1-\sgn M}\right|  \,dy\,.
\end{eqnarray*}
If $M>0 $, then
\begin{eqnarray*}
&  &
\int_0^{z-1}   \left|    \left( \left(1+ z \right)^2-  y ^2\right)^{-  \frac{1}{ 2} }   F \left(  \frac{1}{ 2}+  iM, \frac{1}{ 2}
+  iM;1;\frac{\left(z- 1\right)^2- y ^2}{\left(z+ 1  \right)^2-  y ^2 }\right) \right|  \,dy  
   \leq  
C \quad \mbox{for all}\quad z \in [1,\infty)\,.
\end{eqnarray*}
If $M=0 $, then
\begin{eqnarray*}
&  &
\int_0^{z-1}   \left|    \left( \left(1+ z \right)^2-  y ^2\right)^{-  \frac{1}{ 2} }  \left(1-\ln \left( \frac{4z}{\left(z+ 1  \right)^2-  y ^2 }\right) \right)\right|  \,dy\\
&  \leq &
C +\int_0^{z-1}     \left( \left(1+ z \right)^2-  y ^2\right)^{-  \frac{1}{ 2} }  \left|\ln \left( \frac{4z}{\left(z+ 1  \right)^2-  y ^2 }\right)   \right|  \,dy\,.
\end{eqnarray*}
For $z>R>6$, in the last integral, we split the interval  of integration $[0,z-1]$  into $[0,\sqrt{(z+1)^2-8z}]$ and $[\sqrt{(z+1)^2-8z},z-1] $. 
 For the first interval $y \leq \sqrt{(z+1)^2-8z}$ implies $8z \leq (z + 1)^2 - y^2$. It follows
\begin{eqnarray*} 
\int_0^{\sqrt{(z+1)^2-8z}}   \left( \left(1+ z \right)^2-  y ^2\right)^{-  \frac{1}{ 2} }  \left|\ln \left( \frac{4z}{\left(z+ 1  \right)^2-  y ^2 }\right)   \right|  \,dy 
& \leq &
(\ln z)\int_0^{z-1}     \left( \left(1+ z \right)^2-  y ^2\right)^{-  \frac{1}{ 2} }  \,dy\\
& \leq &
C\ln z\,.
\end{eqnarray*}
For the second interval  $y \in [\sqrt{(z+1)^2-8z},z-1]$ and we obtain   
\begin{eqnarray*} 
\int_{\sqrt{(z+1)^2-8z}}^{z-1}    \left( \left(1+ z \right)^2-  y ^2\right)^{-  \frac{1}{ 2} }  \left|\ln \left( \frac{4z}{\left(z+ 1  \right)^2-  y ^2 }\right)   \right|   \,dy 
& \leq &
C\int_{\sqrt{(z+1)^2-8z}}^{z-1}    \left( \left(1+ z \right)^2-  y ^2\right)^{-  \frac{1}{ 2} }  \,dy\\
& \leq &
C\,.
\end{eqnarray*} 
The lemma is proved. \qed

\subsection{Estimates of  integrals of the kernel functions. Proof of Theorem~\ref{T_large_mass}}
\label{SS3.1}

\begin{lemma}
\label{L3.3}
Let $4m_c^2\geq (3\ell-1)^2   $, then
\begin{eqnarray*}
&  &
\int_0^{ \phi (t)- \phi (b ) }  
 \left| E(r,t;b ; {\cal{M}})\right|  \,dr \leq C|\phi (b  )|^{\frac{1}{1-\ell }}\left(1+\ln   \left(\frac{\phi (b)}{\phi (t)} \right)\right)^{1-\sgn M}
\end{eqnarray*} 
for all  $t \in [1,\infty)$, $1\leq b \leq t$.
\end{lemma}
\medskip

\ndt
{\bf Proof.} 
Denote 
\begin{eqnarray}
\label{toz}
&  &  
z=   \left(\frac{\phi (b)}{\phi (t)} \right)= \left(\frac{ t }{ b }\right)^{\ell-1},\quad y=-\frac{r}{\phi (t)}\,,
\end{eqnarray}
then, from (\ref{Edef}) we derive
\begin{eqnarray*} 
&  &
\int_0^{ \phi (t)- \phi (b ) }  
 \left| E(r,t;b ; {\cal{M}})\right|  \,dr \\ 
& \leq  &
C_\ell\int_0^{ \phi (t)- \phi (b ) }  
 \Bigg|    \phi (b  )^{\frac{ \ell+2i{\cal{M}}}{1-\ell }}
\left( \left(\phi (t)+ \phi (b  )  \right)^2- r^2 \right)^{-  \frac{ i{\cal{M}}}{ 1-\ell }}  F \left(  i \frac{{\cal{M}}}{ 1-\ell },  i\frac{{\cal{M}}}{ 1-\ell };1;\frac{\left( \phi (t)-  \phi (b  )   \right)^2-r^2}{\left( \phi (t)+ \phi (b  ) \right)^2-r^2}\right)\Bigg|  \,dr \\
& \lesssim  &
 |\phi (b  )|^{\frac{1}{1-\ell }}\int_0^{ z- 1 }    
\left( \left(z+ 1  \right)^2- y^2 \right)^{-  \frac{ 1}{ 2}}   \left|  F \left(  \frac{1}{ 2}+  iM , \frac{1}{ 2}+  iM ;1;\frac{\left( z-  1 \right)^2-y^2}{\left( z+ 1 \right)^2-r^2}\right)\right|\,dy\,.
\end{eqnarray*}
Then we apply Lemma~\ref{L3.2} 
and obtain 
\begin{eqnarray*} 
\int_0^{ \phi (t)- \phi (b ) }  
 \left| E(r,t;b ; {\cal{M}})\right|  \,dr  
& \leq  & 
C|\phi (b  )|^{\frac{1}{1-\ell }}(1+\ln (z))^{1-\sgn M}\,.
\end{eqnarray*} 
The lemma is proved. \qed

\begin{lemma}
\label{L3.3b}
Let $4m_c^2\geq (3\ell-1)^2   $, then
\begin{eqnarray*}
&  &
\int_0^{\phi (t)- \phi (1 )}   |K_1 \left(r,t; {\cal{M}} ;1 \right)| \,dr \leq C (1+\ln t)^{1-\sgn M} \quad \mbox{ for all}\quad t \in [1,\infty)\,.
\end{eqnarray*} 
\end{lemma}
\medskip

\ndt
{\bf Proof.} We use notations (\ref{toz}) with $b=1$:
\begin{equation}
\label{zchange}
 z:=\frac{\phi (1)}{\phi (t)}=t^{\ell-1}\in[1,\infty),\quad  y:=-\frac{r}{\phi (t)}
\end{equation}
and obtain
\[ 
\int_0^{\phi (t)- \phi (1 )}   |K_1 \left(r,t; {\cal{M}} ;1 \right)| \,dr   
  \leq  
C \int_0^{z-1}   \left|    \left( \left(1+ z \right)^2-  y ^2\right)^{-  \frac{1}{ 2} }   F \left(  \frac{1}{ 2}+  iM, \frac{1}{ 2}+  iM;1;\frac{\left(z- 1\right)^2- y ^2}{\left(z+ 1  \right)^2-  y ^2 }\right) \right|  \,dy\,.
\]
The application of   Lemma~\ref{L3.2}   completes the proof. \qed

\begin{lemma}
\label{L3.4b}
Let $4m_c^2\geq (3\ell-1)^2   $, then
\begin{eqnarray*}
&  &
\int_0^{\phi (t)- \phi (1 )}   |K_2 \left(r,t; {\cal{M}} ;1 \right)| \,dr \leq C  t^{\frac{\ell-1}{2}}(1+\ln (t))^{1-\sgn M}   \quad \mbox{for all}\quad t \in [1,\infty)\,.
\end{eqnarray*}
\end{lemma}
\medskip

\ndt
{\bf Proof.} 
We use notations  (\ref{zchange}) 
and obtain
\begin{eqnarray*} 
K_2 \left(r,t  ;{\cal{M}};1 \right)
& = &
\frac{2^{2+2 i M} \left(\frac{1}{2}+i M\right) z^{2+2 i M} \left((z+1)^2-y^2\right)^{-\frac{1}{2}-i M}}{(z-1 )^2-y^2}\\
&  &
\times \Bigg\{ (z-1)  F \left(i M+\frac{1}{2},i M+\frac{1}{2};1;\frac{(z-1)^2-y^2}{(z+1)^2-y^2}\right)\\
& &
+\frac{ 2 \left(1-y^2-z^2 \right) }{(z+1)^2-y^2}  F \left(i M+\frac{1}{2},i M+\frac{3}{2};1;\frac{(z-1)^2-y^2}{(z+1)^2-y^2}\right)\Bigg\}\,.
\end{eqnarray*}
By \cite[(28) Sec. 2.8]{B-E},
\begin{eqnarray}
\label{3.20b}
&  &
 F (a+1,a;1;z)=- \frac{(1-a)   F (a-1,a;1;z)+(2 a-1)  F (a,a;1;z)}{a (z-1)} ,
\end{eqnarray}
we derive 
\begin{eqnarray*} 
F \left( \frac{3}{ 2}+ iM,  \frac{1}{ 2}+ iM;1;\frac{\left( z- 1\right) ^2-y^2}{( z+ 1 )^2-y^2}\right) 
&\! = \!&
\frac{\left((z+1)^2-r^2\right)}{2 (1+2 i M) z} \Bigg(2 i M  F \left(i M+\frac{1}{2},i M+\frac{1}{2};1;\frac{(z-1)^2-r^2}{(z+1)^2-r^2}\right)\\
&  &
+\left(\frac{1}{2}-i M\right)  F \left(i M-\frac{1}{2},i M+\frac{1}{2};1;\frac{(z-1)^2-r^2}{(z+1)^2-r^2}\right)\Bigg)\,. 
\end{eqnarray*}
It follows
\begin{eqnarray*}  
K_2 \left(r,t  ;{\cal{M}};1 \right) 
& = &
2^{2+2 i M} \left(\frac{1}{2}+i M\right) \frac{z^{2+2 i M} }{(z-1)^2-y^2}\left((z+1)^2-y^2\right)^{-\frac{1}{2}-i M}\\
&  &
\times \Bigg\{ \left(\frac{2 M \left(1-y^2-z^2 \right)}{(2 M-i) z}+z-1\right) F \left(i M+\frac{1}{2},i M+\frac{1}{2};1;\frac{(z-1)^2-y^2}{(z+1)^2-y^2}\right)\\
& &
+\left(\frac{1}{2}-i M\right)\frac{  1-y^2-z^2    }{(1+2 i M) z}F \left(i M-\frac{1}{2},i M+\frac{1}{2};1;\frac{(z-1)^2-y^2}{(z+1)^2-y^2}\right)\Bigg\} 
\end{eqnarray*}
and
\begin{eqnarray*} 
\int_0^{\phi (t)- \phi (1 )}   |K_2 \left(r,t; {\cal{M}} ;1 \right)| \,dr  
& \lesssim &
 z \int_0^{z- 1}  \frac{1 }{(z-1)^2-y^2}\left((z+1)^2-y^2\right)^{-\frac{1}{2}}\\
&  &
\times \Bigg|  \left(\frac{2 M \left(1-y^2-z^2 \right)}{(2 M-i) z}+z-1\right) F \left(i M+\frac{1}{2},i M+\frac{1}{2};1;\frac{(z-1)^2-y^2}{(z+1)^2-y^2}\right)\\
& &
+\left(\frac{1}{2}-i M\right)\frac{  1-y^2-z^2    }{(1+2 i M) z}F \left(i M-\frac{1}{2},i M+\frac{1}{2};1;\frac{(z-1)^2-y^2}{(z+1)^2-y^2}\right) \Bigg| \,dy \,.
\end{eqnarray*}
For a small positive $\epsilon $ we have
\begin{eqnarray*}
 F \left(\frac{1}{2}+i M,\frac{1}{2}+i M;1;\epsilon \right)
& = & 
1+\left(\frac{1}{2}+i M\right)^2 \epsilon +\frac{1}{64} (2 M-i)^2 (2 M-3 i)^2 \epsilon ^2+O\left(\epsilon ^3\right)\,,\\
 F \left(-\frac{1}{2}+i M,\frac{1}{2}+i M,1,\epsilon\right)
& = &
1-\frac{1}{4} (2 M-i) (2 M+i) \epsilon  
+\frac{1}{64} (2 M-3 i) (2 M-i)^2 (2 M+i) \epsilon ^2+O\left(\epsilon ^3\right)\,.
\end{eqnarray*}
We divide the domain of integration into two zones,
\begin{eqnarray}
\label{Zone1}
Z_1(\varepsilon, z):=\left\{(z, y) \left|  \frac{(z-1)^2-y^2}{(z+1)^2-y^2} \leq \varepsilon, \, 0 \leq y \leq z-1\right\}, \right.\\
\label{Zone2}
Z_2(\varepsilon, z):=\left\{(z, y) \left| \varepsilon \leq \frac{(z-1)^2-y^2}{(z+1)^2-y^2}, \,0 \leq y \leq z-1\right\},\right.
\end{eqnarray} 
and split the integral into two parts,
\begin{eqnarray*}
&  &
\int_0^{\phi (t)- \phi (1 )}   \left|K_2\left(r,t  ;{\cal{M}};1 \right)\right| d r\\
& = &
\int_{(z, y) \in Z_1(\varepsilon, z)}\left|K_2\left(r,t  ;{\cal{M}};1 \right)\right| |\phi (t)|d y+\int_{(z, y) \in Z_2(\varepsilon, z)}\left|K_2\left(r,t  ;{\cal{M}};1 \right)\right||\phi (t)| d y .
 \end{eqnarray*}
In the first zone we have 
\begin{equation}
\label{1zone}
\frac{(z-1)^2-y^2}{(z+1)^2-y^2}\leq \varepsilon <1\,.
\end{equation}
Consequently, the integral over the first zone reads 
\begin{eqnarray*}
&  &
\int_ {Z_1(\epsilon,z)}   |K_2\left(r,t  ;{\cal{M}};1 \right)| \,dr \\ 
& \lesssim &
 z\int_{y \in Z_1(\varepsilon, z)}  \frac{1 }{(z-1)^2-y^2}\left((z+1)^2-y^2\right)^{-\frac{1}{2} }\\
&  &
\times \Bigg|  \left(\frac{2 M \left(1-y^2-z^2 \right)}{(2 M-i) z}+z-1\right)\Bigg\{  1+\left(\frac{1}{2}+i M\right)^2 \epsilon +\frac{1}{64} (2 M-i)^2 (2 M-3 i)^2 \epsilon ^2+O\left(\epsilon ^3\right)\Bigg\} \\
& &
+\left(\frac{1}{2}-i M\right)\frac{  1-y^2-z^2    }{(1+2 i M) z}\Bigg\{ 1-\frac{1}{4} (2 M-i) (2 M+i) \epsilon \\
&  &
+\frac{1}{64} (2 M-3 i) (2 M-i)^2 (2 M+i) \epsilon ^2+O\left(\epsilon ^3\right) \Bigg\}  \Bigg| \,dy   \,. 
 \end{eqnarray*}
If we split the last integral into three parts, then for the first part, which is  defined by
 \begin{eqnarray*} 
A_0 & := & 
  z\int_{y \in Z_1(\varepsilon, z)}  \frac{1 }{(z-1)^2-y^2}\left((z+1)^2-y^2\right)^{-\frac{1}{2} }\\
&  &
\times \Bigg|   \left(\frac{2 M \left(1-y^2-z^2 \right)}{(2 M-i) z}+z-1\right)   
+\left(\frac{1}{2}-i M\right)\frac{  1-y^2-z^2    }{(1+2 i M) z}   \Bigg| \,dy  \,, 
 \end{eqnarray*}
 we have
 \begin{eqnarray*} 
A_0 & \lesssim & 
  \int_0^{z-1}  \frac{1}{  \sqrt{(z+1)^2-y^2}}  \,dy =  \tan ^{-1}\left(\frac{z-1}{ \sqrt{z}}\right) \lesssim 1\quad \mbox{\rm for all}\quad z \in [1,\infty)\,.
 \end{eqnarray*}
 Next, for  $A_1 $, which is defined by
\begin{eqnarray*} 
A_1
& := &
 z\int_{y \in Z_1(\varepsilon, z)}  \frac{1 }{(z-1)^2-y^2}\left((z+1)^2-y^2\right)^{-\frac{1}{2} }\\
&  &
\times \Bigg|  \left(\frac{2 M \left(1-y^2-z^2 \right)}{(2 M-i) z}+z-1\right)\Bigg\{  \left(\frac{1}{2}+i M\right)^2 \epsilon  \Bigg\} \\
& &
+\left(\frac{1}{2}-i M\right)\frac{  1-y^2-z^2    }{(1+2 i M) z}\Bigg\{  -\frac{1}{4} (2 M-i) (2 M+i) \epsilon \Bigg\}  \Bigg| \,dy  \,,
 \end{eqnarray*}
where $\epsilon= \frac{(z-1)^2-y^2}{(z+1)^2-y^2} $,  we obtain
\begin{eqnarray*} 
A_1
& = &
 z\int_{y \in Z_1(\varepsilon, z)} \left((z+1)^2-y^2\right)^{-\frac{3}{2} }\\
&  &
\times \Bigg| \frac{4 M^2 (y-z+1) (y+z-1)-8 i M \left(y^2+z-1\right)+y^2+(z-1) (3 z+1)}{8 z} \Bigg|  \,dy \\
&\lesssim &
 \int_{y \in Z_1(\varepsilon, z)}  \left((z+1)^2-y^2\right)^{-\frac{3}{2} } z^2   \,dy \\ 
&\lesssim &
z^{\frac{1}{2 } }
 \quad \mbox{\rm for all} \quad z \in [1,\infty) \,.
\end{eqnarray*}
Similarly, for $A_2$, which is defined by 
\begin{eqnarray*} 
A_2  
& := &
 z\int_{y \in Z_1(\varepsilon, z)}  \frac{1 }{(z-1)^2-y^2}\left((z+1)^2-y^2\right)^{-\frac{1}{2} }\\
&  &
\times \Bigg|  \left(\frac{2 M \left(1-y^2-z^2 \right)}{(2 M-i) z}+z-1\right)\Bigg\{ \frac{1}{64} (2 M-i)^2 (2 M-3 i)^2 \epsilon ^2\Bigg\} \\
& &
+\left(\frac{1}{2}-i M\right)\frac{  1-y^2-z^2    }{(1+2 i M) z}\Bigg\{ \frac{1}{64} (2 M-3 i) (2 M-i)^2 (2 M+i) \epsilon ^2\Bigg\}  \Bigg| \,dy \,,
 \end{eqnarray*}
we derive
\begin{eqnarray*} 
A_2 
&  \lesssim  &
 z^2\int_{y \in Z_1(\varepsilon, z)}   \left((z+1)^2-y^2\right)^{-\frac{5}{2} }( (z-1)^2-y^2 )\,dy \\
&  \lesssim  &
 z^{\frac{1}{2}  }
\quad \mbox{\rm for all}\quad z \in [1,\infty)\,.
 \end{eqnarray*}
 Finally, for   $A_3$, which is defined by
\begin{eqnarray*} 
A_3
& := &
 z\int_{y \in Z_1(\varepsilon, z)}  \frac{1 }{(z-1)^2-y^2}\left((z+1)^2-y^2\right)^{-\frac{1}{2} }\\
&  &
\times \Bigg|  \left(\frac{2 M \left(1-y^2-z^2 \right)}{(2 M-i) z}+z-1\right)  O\left(\epsilon ^3\right)  
+\left(\frac{1}{2}-i M\right)\frac{  1-y^2-z^2    }{(1+2 i M) z}  O\left(\epsilon ^3\right)    \Bigg| \,dy  \,,
 \end{eqnarray*}
 we obtain 
\begin{eqnarray*} 
A_3
& \lesssim &
 z^2\int_0^{z-1}   \left((z+1)^2-y^2\right)^{-\frac{7}{2} } 
 \left( (z-1)^2-y^2 \right)^2   \,dy   
  =  
 z^2\frac{4 (z-1)^5}{15 \sqrt{z} (z+1)^6}\\
& \lesssim  &
 z^{\frac{1}{2}}\quad \mbox{\rm for all}\quad z \in [1,\infty) \,.
 \end{eqnarray*}
Thus,
\begin{eqnarray*} 
\int_ {Z_1(\epsilon,z)}   |K_2\left(r,t  ;{\cal{M}};1 \right)| \,dr   
& \leq &
C z^{\frac{1}{2}}\quad \mbox{\rm for all}\quad z \in [1,\infty)\,.
\end{eqnarray*}
In the second zone, $Z_2(\varepsilon,z) $, 
we have
\begin{equation}
\label{2zone}
\frac{1}{\left(z- 1\right)^2-y^2 }\leq \frac{1}{\epsilon}\frac{1}{\left(z+ 1\right)^2-y^2}
\end{equation}
and due to (\ref{3.15}) 
\begin{eqnarray*}
\left|  F \left(\frac{1}{ 2}+  iM,\frac{1}{ 2}+ iM;1;\frac{(z-1)^2-y^2}{(z+1)^2-y^2}\right) \right| 
& \leq &  
 C(1+\ln (z))^{1-\sgn M}\,,\\
\left|  F \left(-\frac{1}{ 2}+ iM ,\frac{1}{ 2}+ iM;1;\frac{(z-1)^2-y^2}{(z+1)^2-y^2}\right) \right| 
& \leq & C\,.
 \end{eqnarray*}
Then 
\begin{eqnarray*} 
\int_ {Z_2(\epsilon,z)}   |K_2\left(r,t  ;{\cal{M}};1 \right)|dr 
& \lesssim &
 z \int_{y \in Z_1(\varepsilon, z)}  \frac{1 }{(z-1)^2-y^2}\left((z+1)^2-y^2\right)^{-\frac{1}{2} }\\
&  &
\times \Bigg| \Bigg\{ \left(\frac{2 M \left(1-y^2-z^2 \right)}{(2 M-i) z}+z-1\right) F \left(i M+\frac{1}{2},i M+\frac{1}{2};1;\frac{(z-1)^2-y^2}{(z+1)^2-y^2}\right)\\
& &
+\left(\frac{1}{2}-i M\right)\frac{  1-y^2-z^2    }{(1+2 i M) z}F \left(i M-\frac{1}{2},i M+\frac{1}{2};1;\frac{(z-1)^2-y^2}{(z+1)^2-y^2}\right)\Bigg\} \Bigg| \,dy \\ 
& \lesssim &
 z^2  (1+\ln (z))^{1-\sgn M} \int_0^{z-1}  \frac{1 }{(z-1)^2-y^2}\left((z+1)^2-y^2\right)^{-\frac{1}{2} } 
    \,dy \\   
& \lesssim &
 (1+\ln (z))^{1-\sgn M}z^{\frac{1}{2}} \quad \mbox{\rm for all}\quad z \in [1,\infty)\,.
\end{eqnarray*}
The lemma is proved. \qed

\subsection{Estimates of the solution of the problem without potential and nonlinear terms. Proof of Theorem~\ref{T_large_mass}}
\label{SS3.2}

\begin{proposition}
\label{P3.1}
Let $\ell>1 $ and $4m_c^2\geq (3\ell-1)^2   $ and $\gamma \in \R$. Then the solution $\psi_{ID} $ of the problem 
(\ref{13.12}) for the equation
\begin{eqnarray} 
\label{3.23}
 \frac{\partial^2 \psi }{\partial t^2}
+ \frac{3 \ell}{t}    \frac{\partial \psi }{\partial t}
-  t^{-2\ell}{\mathcal{A}}(x,\partial_x)\psi +  t^{ - 2 }  m_c^2 \psi  = 0\,, 
\end{eqnarray} 
satisfies  
\begin{eqnarray*}
t^\gamma \| \psi_{ID}  (\cdot,t)\|_{ H_{(s)}}   
& \leq &
Ct^{ \gamma - \ell }(\|  \psi_1  \|_{ H_{(s)}} t^{-\frac{ \ell-1}{2}} 
+
  \| \psi _0  \|_{ H_{(s)}} )(1+\ln (t))^{1-\sgn M} \quad \mbox{ for all}\quad t \in [1,\infty)\,. \nonumber
\end{eqnarray*} 
\end{proposition}
\msk

\ndt
{\bf Proof.} 
Indeed,   according to \cite[Theorem 4.1]{JDE2021} if $f=0$, then
\begin{eqnarray*}
 \psi_{ID}  (x,t) 
& = & 
t^{k_-} \int_0^{\phi (t)- \phi (1 )}   K_1 \left(r,t; {\cal{M}} ;1 \right) v_{ \psi_1  (x) }(x,r)\,dr \nonumber \\
&   &
  - k_-t^{k_-} \int_0^{\phi (t)- \phi (1 )}    K_1 \left(r,t; {\cal{M}} ;1 \right) v_{  \psi_0 (x)}(x,r)\,dr \nonumber \\
&   &
+
 t^{k_-} \left(\frac{\phi (t)}{ \phi (1)} \right)^{-   i\frac{{\cal{M}}}{1-\ell} } 
v_{\psi _0}\left(x,\phi (t)- \phi (1) \right)  
+ t^{k_-} \frac{1}{\phi (1)}\int_0^{\phi (t)- \phi (1)}  
K_2\left(r,t  ;{\cal{M}};1 \right)  
v_{\psi _0}(x, r)\,dr  \,. \nonumber
\end{eqnarray*} 
It follows,
\begin{eqnarray*}
\| \psi  (x,t)\|_{ H_{(s)}}  
& \leq &
C t^{\frac{1-3\ell}{2}} \int_0^{\phi (t)- \phi (1 )}  |K_1 \left(r,t; {\cal{M}} ;1 \right) |\|v_{ \psi_1   }(x,r)\|_{ H_{(s)}}\,dr \nonumber \\
&   &
+ c t^{\frac{1-3\ell}{2}} \int_0^{\phi (t)- \phi (1 )} | K_1 \left(r,t; {\cal{M}} ;1 \right)| \|v_{  \psi_0  }(x,r)\|_{ H_{(s)}} \,dr\nonumber \\
&   &
+
C t^{\frac{1-3\ell}{2}}   t  ^{-   \frac{1-\ell}{2} } 
\|v_{\psi _0}\left(x,\phi (t)- \phi (1) \right) \nonumber \|_{ H_{(s)}}\\
&  &
+ Ct^{\frac{1-3\ell}{2}} \int_0^{\phi (t)- \phi (1)}  
|K_2\left(r,t  ;{\cal{M}};1 \right)  |
\|v_{\psi _0}(x, r)  \|_{ H_{(s)}}\,dr\,. \nonumber
\end{eqnarray*} 
Next we apply (\ref{Blplq}) and Lemmas~\ref{L3.3b}, \ref{L3.4b} and obtain
\begin{eqnarray*} 
\| \psi  (x,t)\|_{ H_{(s)}}  
& \leq &
C\|  \psi_1  \|_{ H_{(s)}} t^{\frac{1-3\ell}{2}}(1+\ln (t))^{1-\sgn M} 
+ C \|  \psi_0 \|_{ H_{(s)}}t^{\frac{1-3\ell}{2}}(1+\ln (t))^{1-\sgn M} 
 \nonumber\\
&  &
+ C \| \psi _0  \|_{ H_{(s)}}t^{ - \ell } + C
\| \psi _0  \|_{ H_{(s)}}t^{  - \ell } (1+\ln (t))^{1-\sgn M} \\
& \leq &
C(\|  \psi_1  \|_{ H_{(s)}} t^{-\frac{ \ell-1}{2}} 
+
  \| \psi _0  \|_{ H_{(s)}}) t^{ - \ell } (1+\ln (t))^{1-\sgn M}\,.
\end{eqnarray*}
The proposition is proved. \qed

\subsection{Estimates of the potential term. Proof of Theorem~\ref{T_large_mass}}

\begin{proposition}
\label{P3.2}
If  $\ell>1 $, $\gamma \in \R$, $4m_c^2\geq (3\ell-1)^2   $, and the potential $V$ satisfies condition (\ref{V}), then 
for $M\geq 0$ we have\begin{eqnarray*} 
t^\gamma \| G[V\Phi]\|_{ H_{(s)}}  
& \leq &
 C\epsilon_P \left( \sup_{\tau \in [1,t]} \tau^\gamma \|  \Phi      (x,\tau)\|_{ H_{(s)}} \right) 
\cases{  1 \quad \mbox{if} \quad3\ell-1   >2\gamma \cr
(\ln (t))^{1-\sgn M} t^{\frac{1}{2} ( 2  \gamma  -3 \ell +1)}    \quad \mbox{if} \quad3\ell-1   <2\gamma \cr
 \ln (t) (\ln (t)+1)^{1-\sgn M}  \quad \mbox{if} \quad3\ell-1 =2\gamma   } 
\end{eqnarray*} 
for all $ t \in [1,\infty)$.  
\end{proposition}
\msk

\ndt
{\bf Proof.} Indeed, from Lemma~\ref{L3.3} it follows
\begin{eqnarray}
\label{3.20} 
\| G[F]\|_{ H_{(s)}} \nonumber 
 &  \lesssim  &
  t^{\Re k_-} \int_1 ^{t  } b^{-\Re k_-} \,d b\int_0^{ \phi (t)- \phi (b ) }  
 |E(r,t;b ; {\cal{M}})| \|{\cal E E} [F]   (x,r;b)\|_{ H_{(s)}} \,dr \nonumber \\
 &  \lesssim  &
  t^{\Re k_-} \int_1 ^{t  } b^{-\Re k_-} \| F    (x,b)\|_{ H_{(s)}}\,d b\int_0^{ \phi (t)- \phi (b ) }  
 |E(r,t;b ; {\cal{M}})|  \,dr \nonumber \\
 &  \lesssim &
  t^{\Re k_-} \int_1 ^{t  } b^{-\Re k_-} \| F    (x,b)\|_{ H_{(s)}}|\phi (b  )|^{\frac{1}{1-\ell }}\left(1+\ln   \left(\frac{\phi (b)}{\phi (t)} \right)\right)^{1-\sgn M} \,d b\nonumber \\
 & \lesssim &
   t^{ \frac{1-3\ell }{2}} \int_1 ^{t  } b^{\frac{3\ell+1}{2}} \| F    (x,b)\|_{ H_{(s)}}\left(1+\ln    \left(\frac{t}{b}\right)  \right)^{1-\sgn M} \,d b\,.
\end{eqnarray}
In particular, for the potential term we obtain
\begin{eqnarray*} 
\| G[V\Phi]\|_{ H_{(s)}}  
&  \lesssim  &
  t^{ \frac{1-3\ell }{2}} \int_1 ^{t  } b^{\frac{3\ell+1}{2}-\delta} \|  \Phi  (x,b)\|_{ H_{(s)}}\left(1+\ln    \left(\frac{t}{b}\right)  \right)^{1-\sgn M} \,d b\\
&  \lesssim  &
 \left(\sup_{\tau \in [1,t]} \tau^\gamma \|  \Phi      (x,\tau)\|_{ H_{(s)}}\right) t^{ \frac{1-3\ell }{2}} \int_1 ^{t  } b^{\frac{3\ell+1}{2}-\delta-\gamma} \left(1+\ln    \left(\frac{t}{b}\right)  \right)^{1-\sgn M} \,d b.
\end{eqnarray*}
 If $\delta=2$, then  
\begin{eqnarray*}
t^\gamma \| G[V\Phi]\|_{ H_{(s)}} 
&  \lesssim  &
 \left(\sup_{\tau \in [1,t]} \tau^\gamma \|  \Phi      (x,\tau)\|_{ H_{(s)}}\right) t^{\gamma+ \frac{1-3\ell }{2}} \int_1 ^{t  } b^{\frac{3\ell-3}{2} -\gamma} \left(1+\ln    \left(\frac{t}{b}\right)  \right)^{1-\sgn M} \,d b
\end{eqnarray*}
and with  $ \frac{3\ell-3}{2} -\gamma  \not=-1$ and $M>0$,  we obtain
\begin{eqnarray*} 
 t^{\gamma+ \frac{1-3\ell }{2}} \int_1 ^{t  } b^{\frac{3\ell-3}{2} -\gamma}  \,d b 
&\lesssim &
 \cases{  1 \quad \mbox{if} \quad3\ell-1 -2\gamma >0 ,\cr
 t^{\gamma+ \frac{1-3\ell }{2}}  \quad \mbox{if} \quad3\ell-1 -2\gamma <0\,.}  
\end{eqnarray*}
If $ \frac{3\ell-3}{2} -\gamma  =-1$ and $M>0$, then  
\begin{eqnarray*} 
 t^{\gamma+ \frac{1-3\ell }{2}} \int_1 ^{t  } b^{-1}  \,d b=
t^{\gamma+ \frac{1-3\ell }{2}}  \ln (t) \,. 
\end{eqnarray*}
Thus, if  $M>0$, then
\begin{eqnarray*} 
t^\gamma \| G[V\Phi]\|_{ H_{(s)}}  
& \lesssim &
 \left(\sup_{\tau \in [1,t]} \tau^\gamma \|  \Phi      (x,\tau)\|_{ H_{(s)}}\right) \cases{  1 \quad \mbox{if} \quad3\ell-1 -2\gamma >0 ,\cr
 t^{\gamma+ \frac{1-3\ell }{2}}  \quad \mbox{if} \quad 3\ell-1 -2\gamma <0,\cr
 \ln (t) \quad \mbox{if} \quad3\ell-1 -2\gamma =0  \,.}
\end{eqnarray*}
If $M=0$ and  $ 2 \gamma -3 \ell+1\not= 0$, then we obtain
\begin{eqnarray*}
&  &
t^{\gamma+ \frac{1-3\ell }{2}} \int_1 ^{t  } b^{\frac{3\ell+1}{2}-2-\gamma} \left(1+\ln    \left(\frac{t}{b}\right)  \right)^{1-\sgn M} \,d b\\
& = &
\frac{1}{(2 \gamma -3 \ell+1)^2}\Bigg( -6 \ell t^{\gamma -\frac{3 \ell}{2}+\frac{1}{2}}+4 \gamma  t^{\gamma -\frac{3 \ell}{2}+\frac{1}{2}}-2 t^{\gamma -\frac{3 \ell}{2}+\frac{1}{2}}\\
&  &
-6 \ell \ln (t) t^{\gamma -\frac{3 \ell }{2}+\frac{1}{2}}+4 \gamma  \ln (t) t^{\gamma -\frac{3 \ell}{2}+\frac{1}{2}}+2 \ln (t) t^{\gamma -\frac{3 \ell}{2}+\frac{1}{2}}+6 \ell+2-4 \gamma \Bigg)\,.
\end{eqnarray*}
Then if $M=0$  we obtain 
\[
t^{\gamma+ \frac{1-3\ell }{2}} \int_1 ^{t  } b^{\frac{3\ell+1}{2}-2-\gamma} \left(1+\ln    \left(\frac{t}{b}\right)  \right)^{1-\sgn M} \,d b 
  \lesssim  
\cases{  1 \quad \mbox{if} \quad  2 \gamma -3 \ell+1< 0, \cr
 t^{\gamma -\frac{3 \ell}{2}+\frac{1}{2}} \ln (t)  \quad \mbox{if} \quad  2 \gamma -3 \ell+1> 0,\cr
 \ln (t) (\ln (t)+2)\quad \mbox{if} \quad  2 \gamma -3 \ell+1= 0. }  
\]
Thus, if  $M=0$, then
\begin{eqnarray*} 
t^\gamma \| G[V\Phi]\|_{ H_{(s)}}  
& \lesssim &
 \left(\sup_{\tau \in [1,t]} \tau^\gamma \|  \Phi      (x,\tau)\|_{ H_{(s)}}\right) 
\cases{  1 \quad \mbox{if} \quad3\ell-1 -2\gamma >0, \cr
 t^{\gamma -\frac{3\ell}{2}+\frac{1}{2}} \ln (t)    \quad \mbox{if} \quad3\ell-1 -2\gamma <0,\cr
 \ln (t) (\ln (t)+1)  \quad \mbox{if} \quad3\ell-1 -2\gamma =0. } 
\end{eqnarray*}
Note, if $\gamma \leq \ell$ and $\ell>1$, then $3\ell-1-2\gamma>0$. The lemma is proved. \qed
\medskip

\subsection{Estimates of the nonlinear term. Proof of Theorem~\ref{T_large_mass}}

\begin{theorem}
\label{T3.6}
Assume that   $\ell >1 $, $\gamma \in \R$
 and the semilinear term is Lipschitz continuous with exponent $\alpha \geq 0 $ in the space $H_{(s)} $. 
Then for the operator $G$ with $ M\geq 0 $,   the following estimate holds:
\begin{eqnarray*}
t^\gamma \| G[\Psi(\Phi)]\|_{ H_{(s)}}   
 & \lesssim &
 \left(\sup_{\tau\in [1,t]}  \tau^{\gamma} \| \Phi    (x,\tau)\|_{ H_{(s)}} \right)^{1+\alpha}
 \cases{ t^{\frac{1}{2} ( 2  \gamma  -3 \ell +1)}(\ln (t))^{1-\sgn M}   \quad \mbox{if} \quad  \gamma  > \frac{3 (\ell +1)}{2(\alpha +1)}  ,\cr 
t^{2 - \gamma   \alpha }    \quad \mbox{if} \quad   \gamma  < \frac{3 (\ell +1)}{2 (\alpha +1)},\cr
 t  ^{2 - \gamma   \alpha  }  \ln (t) (\ln (t)+1)^{1-\sgn M}  \quad \mbox{if} \quad \gamma =\frac{3(\ell+1)}{2(1+\alpha)}  }
 \end{eqnarray*}
for all $t \in [1,\infty)$ and all $\Phi \in   H_{(s)}$.
\end{theorem}
\msk

\ndt
{\bf Proof.} Indeed, with $F=\Psi(\Phi)$ from (\ref{3.20}) we derive  
\[
t^\gamma \| G[\Psi(\Phi)]\|_{ H_{(s)}}    
 \lesssim  
 \left(\sup_{\tau\in [1,t]}  \tau^{\gamma} \| \Phi    (x,\tau)\|_{ H_{(s)}} \right)^{1+\alpha}I_{\ell,\alpha,M}(t) \,,
\]
where 
\[
I_{\ell,\alpha,M}(t) :=  t^{\gamma + \frac{1-3\ell }{2}} \int_1 ^{t  } b^{\frac{3\ell+1}{2}- \gamma (1+\alpha) } \left(1+\ln    \left(\frac{t}{b}\right)  \right)^{1-\sgn M} \,d b \,.
\]  
We use the change of variable  $ z=  t / b $, $db=-tz^{-2}dz $, 
and obtain
\[
I_{\ell,\alpha,M}(t) 
  =   
   t  ^{2 - \gamma   \alpha  }\int_{1  }^t   z^{- \frac{3\ell+1}{2}+ \gamma (1+\alpha)-2 } \left(1+\ln   z \right)^{1-\sgn M} \,  dz\,.
\] 
If $M>0 $ and $- \frac{3\ell+1}{2}+ \gamma (1+\alpha)-2\not=-1 $, then 
\[
I_{\ell,\alpha,M}(t) 
  \lesssim    
  \cases{ t^{- \frac{3\ell+1}{2}+ \gamma  +1 } \quad \mbox{if}\quad - \frac{3\ell+1}{2}+ \gamma (1+\alpha)-1>0, \cr
  t  ^{2 - \gamma   \alpha  }\quad \mbox{if}\quad - \frac{3\ell+1}{2}+ \gamma (1+\alpha)-1<0 \,.}   
\] 
If $M>0 $ and $- \frac{3\ell+1}{2}+ \gamma (1+\alpha)-2 =-1 $, then 
$
I_{\ell,\alpha,M}(t) 
  =      t  ^{2 - \gamma   \alpha  }  \ln (t) $.   
Thus, if $M>0 $, then 
\[
t^\gamma \| G[\Psi(\Phi)]\|_{ H_{(s)}}    
  \lesssim  
 \left(\sup_{\tau\in [1,t]}  \tau^{\gamma} \| \Phi    (x,\tau)\|_{ H_{(s)}} \right)^{1+\alpha}
\cases{ t^{- \frac{3\ell+1}{2}+ \gamma  +1 } \quad \mbox{if}\quad - \frac{3\ell+1}{2}+ \gamma (1+\alpha)-1>0, \cr
  t  ^{2 - \gamma   \alpha  }\quad \mbox{if}\quad - \frac{3\ell+1}{2}+ \gamma (1+\alpha)-1<0 ,\cr
t  ^{2 - \gamma   \alpha  }  \ln (t)  \quad \mbox{if}\quad - \frac{3\ell+1}{2}+ \gamma (1+\alpha)-1=0.} 
\]
If $M=0 $ and $- \frac{3\ell+1}{2}+ \gamma (1+\alpha)-2\not=-1 $, then 
\begin{eqnarray*} 
I_{\ell,\alpha,M}(t)  
& = &   
\frac{2}{(2 (\alpha +1) \gamma -3 (\ell +1))^2}  t  ^{2 - \gamma   \alpha  } t^{\frac{1}{2} (-3) (\ell +1)}\\
&  &
\times \Bigg[ \ln (t) (2 (\alpha +1) \gamma -3 (\ell +1)) t^{(\alpha +1) \gamma }+(2 (\alpha +1) \gamma -3 \ell -5) \left(t^{(\alpha +1) \gamma }-t^{\frac{3 (\ell +1)}{2}}\right)\Bigg]\\
& \lesssim &  
 t^{\frac{1}{2} (-2 \alpha  \gamma -3 \ell +1)} \Bigg[ \ln (t)  t^{(\alpha +1) \gamma }+|2 (\alpha +1) \gamma -3\ell -5|  t^{\frac{3 (\ell +1)}{2}} \Bigg]\\
& \lesssim &  
 t^{\frac{1}{2} (-2 \alpha  \gamma -3 \ell +1)}
 \cases{   t^{(\alpha +1) \gamma }\ln (t) \quad \mbox{if} \quad (\alpha +1) \gamma  \geq \frac{3 (\ell +1)}{2} ,\cr 
t^{\frac{3 (\ell +1)}{2}}   \quad \mbox{if} \quad   (\alpha +1) \gamma  < \frac{3 (\ell +1)}{2} \,.}
\end{eqnarray*} 
If $M=0 $ and $- \frac{3\ell+1}{2}+ \gamma (1+\alpha)-2 =-1 $, then 
\[ 
I_{\ell,\alpha,M}(t)  
  =     
   t  ^{2 - \gamma   \alpha  }  \frac{1}{2} \ln (t) (\ln (t)+2)\lesssim   t  ^{2 - \gamma   \alpha  }  \ln (t) (\ln (t)+1)\,.
\] 
Thus, if $M=0 $, then
\begin{eqnarray*} 
&  &
t^\gamma \| G[\Psi(\Phi)]\|_{ H_{(s)}}  \\
 & \lesssim &
 \left(\sup_{\tau\in [1,t]}  \tau^{\gamma} \| \Phi    (x,\tau)\|_{ H_{(s)}} \right)^{1+\alpha}
  \cases{ t^{\frac{1}{2} ( 2  \gamma  -3 \ell +1)}\ln (t)   \quad \mbox{if} \quad (\alpha +1) \gamma  > \frac{3 (\ell +1)}{2}  ,\cr 
t^{2 - \gamma   \alpha }    \quad \mbox{if} \quad   (\alpha +1) \gamma  < \frac{3 (\ell +1)}{2},\cr
 t  ^{2 - \gamma   \alpha  }  \ln (t) (\ln (t)+1)  \quad \mbox{if} \quad - \frac{3\ell+1}{2}+ \gamma (1+\alpha)-1 =0  \,.}
 \end{eqnarray*}
The lemma is proved. \qed

\subsection{Proof of global existence. Proof of Theorem~\ref{T_large_mass}}
\label{SS3.5}

We intend to apply Banach's fixed-point theorem. We use the Lipschitz condition ($\mathcal L$) to estimate nonlinear
term. First, we consider the following integral equation
\begin{equation} 
\label{5.1}
\psi  (x,t)
 = 
\psi  _{ID}(x,t) + G[V  \psi ] (x,t)+
G[ F(\cdot)\Psi(\cdot ,\psi ) ] (x,t)    \,, 
\end{equation}  
where the function $ \psi _{ID} (x,t) \in C([1,\infty);H_{(s)})$ is given.  We set $ a_{acc}=1$ and $c=1$.
\medskip

We consider the Cauchy problem in the Sobolev space $ H_{(s)}$ with $ s > 3/2$, which is an algebra.

Denote $\widetilde{C}^\ell ([1,T]; H_{(s)}) $ the complete subspace of $C^\ell([1,T]; H_{(s)})$ of all   functions $f=f(x,t)$ with supp\,$f \subset \{(x,t) \in \R^3\times[1,\infty)\,|\,|x|\geq R_{ID}-( \phi (t)-\phi(1)) \,\}$.  
  For every $T>1$, due to Section~\ref{S2}, the operator $G$ maps
\[
G \,:\, \widetilde{C}([1,T]; H_{(s)})\longrightarrow \widetilde{C}^2([1,T]; H_{(s)})
\]
 continuously.   
Consider the mapping
\begin{equation}
\label{14.5}
S [\Phi] 
= \psi_{ID}+G[V\Phi] +G[F\Psi(\Phi)]\,,
\end{equation}
where  
$
\psi_{ID} \in X({R,H_{(s)},\gamma})$. 
The function $
\psi_{ID}  
$ 
can be  generated by initial data      (\ref{13.13}),(\ref{13.12}). In that case it  belongs to $X(R,H_{(s)},\gamma) $ if $\gamma \leq \ell$ and 
\[
t^{\gamma }\| \psi_{ID}(x,t )\|_{H_{(s)}}
  \leq  
C ( \| \psi _0 \|_{H_{(s)}}+ \| \psi _1 \|_{H_{(s)}})\,.
\]
It is easily seen that   supp\,$S [\Phi] \subset$ supp\,$ \Phi $ if   supp\,$ \Phi  \subset$ supp\,$ \psi_{ID} $. We claim that since  $\ell >1 $ and $ \gamma \leq \ell$, if $ \Phi \in X({R,H_{(s)},\gamma})$ with  $ \gamma\leq\frac{3\ell+1}{ \alpha } $ and if  supp~$\Phi \subseteq \{(x,t) \in \R^3\times[0,\infty)\,|\,|x|> R_{ID}-( \phi (t)-\phi(1)) \,\} $, then 
$S[\Phi] \in   X({R,H_{(s)},\gamma})$.  Moreover, $S$ 
is a contraction, provided that   $\varepsilon  $, $\varepsilon_P  $, and $R$ are sufficiently small. 

First, we note that due to Theorem~\ref{T3.6} and Propositions~\ref{P3.1}, \ref{P3.2}, for the values of $\alpha$, $\ell$,  and $\gamma$ assumed in the theorem, we have $\gamma-\ell \leq 0 $ and, consequently,
\begin{eqnarray*} 
t^{\gamma }\|S [\Phi] \|_{H_{(s)}}
& \leq &
 t^{\gamma }\| \psi_{ID}\|_{H_{(s)}}+t^{\gamma }\|G[V\Phi]\|_{H_{(s)}} +t^{\gamma }\|G[F\Psi(\Phi)]\|_{H_{(s)}}\\
& \leq &
Ct^{ \gamma - \ell }(\|  \psi_1  \|_{ H_{(s)}} t^{-\frac{ \ell-1}{2}} 
+
  \| \psi _0  \|_{ H_{(s)}} )(1+\ln (t))^{1-\sgn M} \\
  &  &
  +C\epsilon_P \left( \sup_{\tau \in [1,t]} \tau^\gamma \|  \Phi      (x,\tau)\|_{ H_{(s)}} \right) 
\cases{  1 \quad \mbox{if} \quad3\ell-1   >2\gamma \cr
(\ln (t))^{1-\sgn M} t^{\frac{1}{2} ( 2  \gamma  -3 \ell +1)}     \quad \mbox{if} \quad3\ell-1   <2\gamma \cr
 \ln (t) (\ln (t)+1)^{1-\sgn M}  \quad \mbox{if} \quad3\ell-1 =2\gamma   } \\
 &  &
 +
 C_{NT}\left(\sup_{\tau\in [1,t]}  \tau^{\gamma} \| \Phi    (x,\tau)\|_{ H_{(s)}} \right)^{1+\alpha}
 \cases{ t^{\frac{1}{2} ( 2  \gamma  -3 \ell +1)}(\ln (t))^{1-\sgn M}   \quad \mbox{if} \quad  \gamma  > \frac{3 (\ell +1)}{2(\alpha +1)}  ,\cr 
t^{2 - \gamma   \alpha }    \quad \mbox{if} \quad   \gamma  < \frac{3 (\ell +1)}{2 (\alpha +1)},\cr
 t  ^{2 - \gamma   \alpha  }  \ln (t) (\ln (t)+1)^{1-\sgn M}  \quad \mbox{if} \quad \gamma =\frac{3(\ell+1)}{2(1+\alpha)} \,.}
\end{eqnarray*}

Hence, to get a contraction we set all time-depending factors to be bounded for all times.  {For $M>0$:} we choose  
  { $\frac{3 (\ell +1)}{2(\alpha +1)} < \gamma \leq \ell \Longrightarrow$ $\alpha >\frac{\ell +3}{2\ell} $}; 
 { $ \frac{2}{\alpha}\leq \gamma    < \frac{3 (\ell +1)}{2(\alpha +1)}  $, $\gamma \leq \ell  $, $  \alpha>\frac{4}{3\ell-1}$}; 
 {$\gamma =\frac{3(\ell+1)}{2(1+\alpha)} $,  $ \alpha >\frac{4}{3\ell-1}.  $} 
  {For $M=0$  we choose:} 
  { $\frac{3 (\ell +1)}{2(\alpha +1)} < \gamma < \ell \Longrightarrow $ $\alpha >\frac{\ell +3}{2\ell} $}; 
 { $ \frac{2}{\alpha}\leq \gamma    < \frac{3 (\ell +1)}{2(\alpha +1)}  $, $\gamma < \ell  $, $  \alpha>\frac{4}{3\ell-1}$} ; 
 {  $\gamma =\frac{3(\ell+1)}{2(1+\alpha)} $,  $\alpha >\frac{4}{3\ell-1}, \gamma < \ell  $}.  
Note that $\frac{\ell +3}{2\ell}> \frac{4}{3\ell-1}$.

Thus, we choose $\gamma$ in accordance to the cases above, and the last inequality proves that the operator $S$ maps $X({R,s,\gamma})$ into itself if 
\[
\| \psi _0 \|_{H_{(s)}}+ \| \psi _1 \|_{H_{(s)}} \leq \varepsilon \,,
\]
$\varepsilon_P  $, $\varepsilon  $, and $R$ are sufficiently small, namely, if $  \varepsilon_P C_{pt}<1 $ and 
\[
\frac{1}{ 1- \varepsilon_P C_{pt} } C_{id} 
\varepsilon  + \frac{1}{ 1- \varepsilon_P C_{pt}  }C_{nt} R^{\alpha +1} < R  .
\]
We derive the contraction property of $S$  from  
\begin{eqnarray*} 
d( S[\Phi_1 ],S[\Phi_2  ] ) 
& \leq &
 d( \Phi_1,\Phi_2 )( C_{pt} \varepsilon_0  + C_{nt} R^\alpha )   \,.
\end{eqnarray*}
Indeed,  due to Theorem~\ref{T3.6}, Proposition~\ref{P3.2}, and conditions of theorem we have
\begin{eqnarray*} 
&  &
t^{\gamma }\| S[\Phi_1 ](\cdot , t) -  S[\Phi_2 ](\cdot , t) \|_{H_{(s)}} \\
& \leq &
t^{\gamma }\|G[V(\Phi_1-\Phi_2)\|_{H_{(s)}} +t^{\gamma }\|G[F(\Psi(\Phi_1)-\Psi(\Phi_2))]\|_{H_{(s)}}\\
& \leq &
C_{pt} \varepsilon_P   \sup_{\tau \in [1,t]} \tau^{\gamma} \| ( \Phi_1  - \Phi_2 )    (x,\tau)\|_{ H_{(s)}}
+Ct^{\gamma }\|G[F(\Psi(\Phi_1)-\Psi(\Phi_2))]\|_{H_{(s)}}  ,
\end{eqnarray*}
while according to above used arguments  and  due to the condition ($\mathcal L$), we obtain
\begin{eqnarray*}
&  &
\| G[(\Psi(\Phi_1)-\Psi(\Phi_2))]\|_{ H_{(s)}}\\
 & \leq &
 Ct^{ \Re k_-} \int_1 ^{t  } b^{-\Re k_-} \| (\Psi(\Phi_1)-\Psi(\Phi_2))  (x,b)\|_{ H_{(s)}} |\phi (b) |^{\frac{1}{1-\ell} } \left( 1+ \ln \frac{\phi (b) }{\phi (t)}\right)^{ 1- \sgn M}  \,d b  \nonumber \\
& \leq &
 Ct^{\frac{1-3\ell}{2} } \int_1 ^{t  } b^{\frac{1+3\ell}{2} }  \|  \Phi_1 (x,b)- \Phi_2   (x,b)\|_{ H_{(s)}}\left( \|  \Phi_1 (x,b)\|_{ H_{(s)}}^\alpha+\| \Phi_2   (x,b)\|_{ H_{(s)}}^\alpha\right)   \left( 1+ \ln \frac{\phi (b) }{\phi (t)}\right)^{ 1- \sgn M} \,d b  \\
& \leq &
  Ct^{\frac{1-3\ell}{2} } \int_1 ^{t  } b^{\frac{1+3\ell}{2} -\gamma(1+\alpha)}\left( \sup_{b \in [1,t]} b^\gamma \|  \Phi_1 (x,b)- \Phi_2   (x,b)\|_{ H_{(s)}}\right) \\
&  &
\times \left( \left( \sup_{b \in [1,t]} b^\gamma \|  \Phi_1 (x,b)\|_{ H_{(s)}}\right) ^\alpha+\left( \sup_{b \in [1,t]} b^\gamma \| \Phi_2   (x,b)\|_{ H_{(s)}}\right) ^\alpha\right)\left( 1+ \ln \frac{\phi (b) }{\phi (t)}\right)^{ 1- \sgn M}\,d b \,.
 \end{eqnarray*}
 Further, if 
 \[
 d( \Phi_1,\Phi_2 ) :=  \sup_{b \in [1,\infty]} b^\gamma \|  \Phi_1 (x,b)- \Phi_2   (x,b)\|_{ H_{(s)}} 
 \]
 and
 \[
 R:=\max\left\{    \sup_{b \in [1,t]} b^\gamma \|  \Phi_1 (x,b)\|_{ H_{(s)}}  \, ,  \, \sup_{b \in [1,t]} b^\gamma \| \Phi_2   (x,b)\|_{ H_{(s)}}  \right\}\,,
 \]
 then
\[
\| G[(\Psi(\Phi_1)-\Psi(\Phi_2))]\|_{ H_{(s)}}
  \leq  
  Cd( \Phi_1,\Phi_2 )R^\alpha  t^{\frac{1-3\ell}{2} } \int_1 ^{t  } b^{\frac{1+3\ell}{2} -\gamma(1+\alpha)}\left( 1+ \ln \frac{\phi (b) }{\phi (t)}\right)^{ 1- \sgn M}\,d b 
 \]
 and, consequently, we obtain
 \[ 
t^\gamma \| G[(\Psi(\Phi_1)-\Psi(\Phi_2))]\|_{ H_{(s)}} 
  \leq  
  Cd( \Phi_1,\Phi_2 )R^\alpha \cases{ t^{\frac{1}{2} ( 2  \gamma  -3 \ell +1)}(\ln (t))^{1-\sgn M}   \quad \mbox{if} \quad  \gamma  > \frac{3 (\ell +1)}{2(\alpha +1)}  ,\cr 
t^{2 - \gamma   \alpha }    \quad \mbox{if} \quad   \gamma  < \frac{3 (\ell +1)}{2 (\alpha +1)},\cr
 t  ^{2 - \gamma   \alpha  }  \ln (t) (\ln (t)+1)^{1-\sgn M}  \quad \mbox{if} \quad \gamma =\frac{3(\ell+1)}{2(1+\alpha)} \,.} 
 \] 
Thus, in all cases   we have
\[ 
t^{\gamma }\| S[\Phi_1 ](\cdot , t) -  S[\Phi_2 ](\cdot , t) \|_{H_{(s)}} 
  \leq  
t^{\gamma }\|G[V(\Phi_1-\Phi_2)\|_{H_{(s)}} +t^{\gamma }\|G[F(\Psi(\Phi_1)-\Psi(\Phi_2))]\|_{H_{(s)}} \,.  
\]
 Finally,
\[ 
t^{\gamma }\| S[\Phi_1 ](\cdot , t) -  S[\Phi_2 ](\cdot , t) \|_{H_{(s)}} 
  \leq  
C_{pt} \varepsilon_P   d( \Phi_1,\Phi_2 )+ C_{nt}d( \Phi_1,\Phi_2 )R^\alpha   \,.
\]
Then we choose $\| \psi _0 \|_{H_{(s)}}+ \| \psi _1 \|_{H_{(s)}} \leq \varepsilon$ with sufficiently small $\varepsilon $  and $R$ such that $C_{pt} \varepsilon_P +C_{nt} R^\alpha<1 $. Banach's fixed point theorem completes the proof. \qed

\section{Case of  small mass. Proof of Theorem~\ref{T_small_mass} }
\label{S4}

Since $\ell >1 $, for the case of real small mass,      $m_c \leq [0,(3\ell-1)/2) $. 
We choose $k_+$ in the Liouville substitution:  
\begin{eqnarray*}
&  &
\psi  (x,t)=  t^{k_+} u (x,t)=t^{\frac{1}{2} \left(1-3 \ell  +\sqrt{(1-3\ell  )^2 -4 m_c^2}\right)} u (x,t), \\
&  &
\psi_0  (x)=   u_0 (x ),\quad \psi_1  (x )=  u_1(x )+ k_+ u_0 (x ), \quad    
u_1(x )=  \psi_1  (x)- k_+ \psi_0 (x), \quad
g(x,t)= t^{-k_+}f(x,t) .
\end{eqnarray*}
  Hence, by Theorem~\ref{T2.4}, 
\begin{eqnarray*}
 \psi  (x,t) 
& = &
 t^{k_+}  \int_1 ^{t  }  b^{-k_+}\,d b\int_0^{ \phi (t)- \phi (b ) }  
 E(r,t;b ; {\cal{M}}) v_{ f}  \left(x,r;b  \right) \,dr  \\
&   &
+ t^{k_+} \int_0^{\phi (t)- \phi (1 )}    K_1 \left(r,t; {\cal{M}} ;1 \right) v_{ \psi_1  (x) }(x,r)\,dr \nonumber \\
&   &
  - k_+t^{k_+} \int_0^{\phi (t)- \phi (1 )}   K_1 \left(r,t; {\cal{M}} ;1 \right) v_{  \psi_0 (x)}(x,r)\,dr \nonumber \\
&   &
+
 t^{k_+} \left(\frac{\phi (t)}{ \phi (1)} \right)^{-   i\frac{{\cal{M}}}{1-\ell} } 
v_{\psi _0}\left(x,\phi (t)- \phi (1) \right) 
+ t^{k_+} \frac{1}{\phi (1)}\int_0^{\phi (t)- \phi (1)}  
K_2\left(r,t  ;{\cal{M}};1 \right)  
v_{\psi _0}(x, r)\,dr  \,. \nonumber
\end{eqnarray*} 
Since 
$  -\frac{\ell}{\ell-1}<  i \frac{{\cal{M}_+}}{ 1-\ell }  
 <\frac{1}{2} $,  
by use of (\ref{Blplq}) we derive
\begin{eqnarray*}
\|  u(x,t )\|_{ H_{(s)}}
& \leq &
C_s t^{-\frac{1}{  2  }\left(1-\frac{\sqrt{(1-3 \ell )^2-4 m_c^2}}{\ell-1}  \right)}  
\| u _0 \|_{ H_{(s)}}
+C_s\|  u _1  \|_{ H_{(s)}}\int_0^{\phi (t)- \phi (1 )}   |K_1 \left(r,t; {\cal{M}} ;1 \right)| \,dr
\\
&  &
+C_s\| u _0 \|_{ H_{(s)}}\int_0^{\phi (t)- \phi (1)} \left|K_2\left(r,t  ;{\cal{M}};1 \right)\right|  
\,dr \,.
\end{eqnarray*} 
For the real small mass, the condition  $1-\frac{\sqrt{(1-3 \ell )^2-4 m_c^2}}{\ell-1}\geq 0$ is equivalent to $4 m_c^2 \geq 4\ell(2\ell-1)$. 
Thus,
\begin{eqnarray*} 
\|  u(x,t )\|_{ H_{(s)}} 
& \leq &
C_s \| u _0 \|_{ H_{(s)}} \cases{  t^{-\frac{1}{  2  }\left(1-\frac{\sqrt{(1-3 \ell )^2-4 m_c^2}}{\ell-1}  \right)}
\quad \mbox{if}\quad 4 m_c^2 <  4\ell(2\ell-1),\cr
1 \quad \mbox{if}\quad 4 m_c^2 \geq 4\ell(2\ell-1) , }  
\\
&   &
+C_s\|  u _1  \|_{ H_{(s)}}\int_0^{\phi (t)- \phi (1 )}   |K_1 \left(r,t; {\cal{M}} ;1 \right)|  dr 
+C_s\| u _0 \|_{ H_{(s)}}\int_0^{\phi (t)- \phi (1)} \left|K_2\left(r,t  ;{\cal{M}};1 \right)\right|  
dr .
\end{eqnarray*}

\subsection{Estimates of the kernel functions. Proof of Theorem~\ref{T_small_mass}}
\label{SS4.2}

\begin{lemma}
\label{L4.1} 
For $4 m_c^2 < (1-3 \ell )^2 $ and the given above parameters ${\cal{M}} $ and $\ell $, 
\begin{eqnarray*}
 \hspace{-1.7cm}\left| F \left(  \frac{1}{2},     i \frac{{\cal{M}}}{1-\ell },\frac{3}{2};\frac{\left(z-1\right)^2 }{\left(z+ 1   \right)^2  }\right)\right|  
\! & \! \leq \! & \! 
C \,,\\  
\hspace{-1.7cm} \left|F \left(\frac{1}{2},1+ i\frac{{\cal{M}}}{1-\ell };\frac{3}{2};\frac{(z-1)^2}{(z+1)^2}\right) \right| 
 \! & \! \leq \! & \! 
C \cases{1 \quad \mbox{if}\quad 4m_c^2 \leq 4\ell(2\ell-1) \,,\cr
 z^{ \frac{1}{  2}-\frac{1}{  2(\ell - 1 ) } \sqrt{(1-3 \ell )^2-4 m_c^2} )}  \,\,
 \mbox{if } \,\,  4m_c^2 > 4\ell(2\ell-1)}  
 \end{eqnarray*}
for all $z \in[1,\infty) $.
\end{lemma}
\medskip

\ndt
{\bf Proof.} The first inequality follows from   \cite[(14) Sec.2.1]{B-E}   since $1 +  i \frac{{\cal{M}}}{\ell-1} > 0 $.
Similarly, for the second one we use   \cite[(14),(23) Sec.2.1]{B-E} and derive
 \begin{eqnarray*} 
\left|F \left(\frac{1}{2},1+ i\frac{{\cal{M}}}{1-\ell };\frac{3}{2};\frac{(z-1)^2}{(z+1)^2}\right) \right| 
 &  \leq &
 C\left(1+ \left(\frac{(z+1)^2}{4z} \right)^{  i\frac{{\cal{M}}}{1-\ell } }\right)    
 \end{eqnarray*} 
  for all $z \in[1,\infty)$. The lemma is proved. \qed

\begin{lemma}
\label{L4.2E}
If $\ell>1$ and  $4 m_c^2 < (1-3 \ell )^2 $,  then for the given above parameters ${\cal{M}} $ and $\ell $, 
\begin{eqnarray*} 
\int_0^{\phi (t)- \phi (b )}   |E \left(r,t;b; {\cal{M}}   \right)| \,dr  
& \lesssim  &
 b ^{1+ \sqrt{(1-3 \ell )^2-4 m_c^2} }  
   | \phi (b)+\phi (t) |^{ \frac{\sqrt{(1-3 \ell )^2-4 m_c^2}}{ \ell-1 }}  
\end{eqnarray*} 
  for all $t \in[1,\infty)$.
\end{lemma}
\medskip

\ndt
{\bf Proof.} From (\ref{Edef})
\begin{eqnarray*}
\int_0^{\phi (t)- \phi (b )}   |E \left(r,t;b; {\cal{M}}   \right)| \,dr 
& \lesssim  &
\int_0^{\phi (t)- \phi (b)}   \Bigg|   \phi (b )^{\frac{ \ell+2i{\cal{M}}}{1-\ell }}
\left( \left(\phi (t)+ \phi (b  )  \right)^2- r^2 \right)^{-  \frac{ i{\cal{M}}}{ 1-\ell }} \\
&  &
\times  F \left(  i \frac{{\cal{M}}}{ 1-\ell },  i\frac{{\cal{M}}}{ 1-\ell };1;\frac{\left( \phi (t)-  \phi (b )   \right)^2-r^2}{\left( \phi (t)+ \phi (b ) \right)^2-r^2}\right)\Bigg| \,dr  .
\end{eqnarray*}
In the notations  (\ref{toz}), we write  
\begin{eqnarray*}
&  &
\int_0^{\phi (t)- \phi (b )}   |E \left(r,t;b; {\cal{M}}   \right)| \,dr  \\
& \lesssim  &
 | \phi (b )^{\frac{ \ell+2i{\cal{M}}}{1-\ell }} \phi (t)^{1-  2\frac{ i{\cal{M}}}{ 1-\ell }}| \int_0^{z- 1}   \Bigg| 
\left( \left(z+ 1  \right)^2- y^2 \right)^{-  \frac{ i{\cal{M}}}{ 1-\ell }} F \left(  i \frac{{\cal{M}}}{ 1-\ell },  i\frac{{\cal{M}}}{ 1-\ell };1;\frac{\left( z-  1   \right)^2-y^2}{\left(z+ 1 \right)^2-y^2}\right)\Bigg| \,dy  \\
\end{eqnarray*} 
and then,  since $ 2i \frac{{\cal{M}}}{1-\ell}<1 $, $\frac{ \ell+2i{\cal{M}}}{1-\ell }= -\frac{1}{\ell-1}-\frac{\sqrt{(1-3 \ell )^2-4 m_c^2}}{\ell-1}$, and $1-  2\frac{ i{\cal{M}}}{ 1-\ell }= \frac{\sqrt{(1-3 \ell )^2-4 m_c^2}}{\ell-1}$, we obtain
\[ 
\int_0^{\phi (t)- \phi (b )}   |E \left(r,t;b; {\cal{M}}   \right)| \,dr   
  \lesssim  
 \left| \phi (b )^{-\frac{1}{\ell-1}-\frac{\sqrt{(1-3 \ell )^2-4 m_c^2}}{\ell-1}} \phi (t)^{\frac{\sqrt{(1-3 \ell )^2-4 m_c^2}}{\ell-1}}\right|  \int_0^{z- 1}    
\left( \left(z+ 1  \right)^2- y^2 \right)^{-  \frac{ i{\cal{M}}}{ 1-\ell }}   \,dy  . 
\]
On the other hand, since $\frac{1}{2}+ \frac{ i{\cal{M}}}{ 1-\ell }=1- \frac{\sqrt{(1-3 \ell )^2-4 m_c^2}}{2(\ell-1)}<\frac{3}{2}$, we have
\begin{eqnarray*}   
  \int_0^{z- 1}    
\left( \left(z+ 1  \right)^2- y^2 \right)^{-  \frac{ i{\cal{M}}}{ 1-\ell }}   \,dy     
& = &
(z-1) (z+1)^{-  2 \frac{ i{\cal{M}}}{ 1-\ell }}  F \left(\frac{1}{2},\frac{ i{\cal{M}}}{ 1-\ell };\frac{3}{2};\frac{(z-1)^2}{(z+1)^2}\right)\\ 
& \lesssim &
(z-1) (z+1)^{-1+ \frac{\sqrt{(1-3 \ell )^2-4 m_c^2}}{ \ell-1 }} \,.
\end{eqnarray*}
It follows
\begin{eqnarray*} 
\int_0^{\phi (t)- \phi (b )}   |E \left(r,t;b; {\cal{M}}   \right)| \,dr   
& \lesssim  &
  | \phi (b )|^{-\frac{1}{\ell-1}-\frac{\sqrt{(1-3 \ell )^2-4 m_c^2}}{\ell-1}}  
   | \phi (b)+\phi (t)|^{ \frac{\sqrt{(1-3 \ell )^2-4 m_c^2}}{ \ell-1 }}\,.  
\end{eqnarray*} 
 The lemma is proved. \qed

\begin{lemma}
If $\ell>1$ and  $4 m_c^2 < (1-3 \ell )^2 $,  then for the given above parameters ${\cal{M}} $ and $\ell $ 
\begin{eqnarray*}
&  &
\int_0^{\phi (t)- \phi (1 )}   |K_1 \left(r,t; {\cal{M}} ;1 \right)| \,dr \leq C \quad \mbox{for all}\quad t \in[1,\infty)\,.
\end{eqnarray*} 
\end{lemma}
\medskip

\ndt
{\bf Proof.} According to (\ref{K1def}) and   $ 2i \frac{{\cal{M}}}{1-\ell}<1 $ we obtain
\begin{eqnarray*}
&  &
\int_0^{\phi (t)- \phi (1 )}   |K_1 \left(r,t; {\cal{M}} ;1 \right)| \,dr \\
& \leq &
C \int_0^{\phi (t)- \phi (1 )}    \left|    \left( \left(\phi (t)+ \phi (1 )  \right)^2-  r ^2\right)^{-  i \frac{{\cal{M}}}{1-\ell}}   F \left(  i \frac{{\cal{M}}}{1-\ell},  i \frac{{\cal{M}}}{1-\ell};1;\frac{\left(\phi (t)- \phi (1 )\right)^2- r ^2}{\left(\phi (t)+ \phi (1 )   \right)^2-  r ^2 }\right) \right| \,dr\\
& \leq &
C \int_0^{\phi (t)- \phi (1 )}        \left( \left(\phi (t)+ \phi (1 )  \right)^2-  r ^2\right)^{   i \frac{{\cal{M}}}{\ell-1}}   \,dr
\,.
\end{eqnarray*}
Here  $\phi (t)- \phi (1)\geq0$. 
In the notations (\ref{zchange}),  
  since $ 2i \frac{{\cal{M}}}{1-\ell}<1 $, we obtain 
\begin{eqnarray*} 
\int_0^{\phi (t)- \phi (1 )}   |K_1 \left(r,t; {\cal{M}} ;1 \right)| \,dr  
& \lesssim  &
  |(\phi (t))^{1+  2 i \frac{{\cal{M}}}{\ell-1}}| 
\int_0^{z- 1}       \left( \left(z+1  \right)^2-  y ^2\right)^{   i \frac{{\cal{M}}}{\ell-1}}   \,dy \,.
\end{eqnarray*}
Then according to \cite[(14) Sec. 2.1]{B-E} and since $1 +  i \frac{{\cal{M}}}{\ell-1} > 0 $, we have
\begin{eqnarray*} 
\int_0^{\phi (t)- \phi (1 )}   |K_1 \left(r,t; {\cal{M}} ;1 \right)| \,dr 
&\lesssim  & 
 |(\phi (t))^{1+  2 i \frac{{\cal{M}}}{\ell-1}}| (z-1) (z+1)^{   2i \frac{{\cal{M}}}{\ell-1}} 
F \left(  \frac{1}{2}, i\frac{{\cal{M}}}{1-\ell },\frac{3}{2};\frac{\left(z-1\right)^2 }{\left(z+ 1   \right)^2  }\right)\\
& \lesssim & 
1 \,.
\end{eqnarray*} 
The lemma is proved. \qed

\begin{lemma}If  $ \ell>1$ and $4 m_c^2<(1-3 \ell )^2$, then  we have
\begin{eqnarray*}
&  &
\int_0^{\phi (t)- \phi (1 )}   |K_2\left(r,t  ;{\cal{M}};1 \right)| \,dr \leq 
C\cases{ \dsp t  ^{ \frac{\ell - 1}{  2}-\frac{1}{  2} \sqrt{(1-3 \ell )^2-4 m_c^2} )}  \quad   \mbox{if }\quad 4m_c^2>  4\ell(2\ell-1),   \cr
1 \quad   \mbox{if }\quad 4m_c^2\leq 4\ell(2\ell-1)  }
\end{eqnarray*} 
for all $t \in[1,\infty)$.
\end{lemma}
\medskip

\ndt
{\bf Proof.} According to the definition of the kernel $K_2$ we have
\begin{eqnarray*} 
\int_0^{\phi (t)- \phi (1 )}   |K_2\left(r,t  ;{\cal{M}};1 \right)| \,dr   
& \leq &
\int_0^{\phi (t)- \phi (1 )}   \Bigg|- i  {\cal{M}} { 2^{1+2 i  \frac{{\cal{M}}}{1-\ell }}}  \phi^{  \frac{2i{\cal{M}}}{1-\ell }} (1 )
\left(\left(\phi (t)+ \phi (1 )\right)^2-r^2\right)^{-  i\frac{{\cal{M}}}{1-\ell }}\\
&  &
\times \Bigg[  \frac{\phi^2 (1 )  -\phi (t)\phi (1 )    }{ r^2-\left(\phi (t)- \phi (1)\right)^2 }  
F \left(  i \frac{{\cal{M}}}{1-\ell },  i \frac{{\cal{M}}}{1-\ell };1;\frac{\left( \phi (t)- \phi (1) \right) ^2-r^2}{(  \phi (t)+ \phi (1 ) )^2-r^2}\right)\nonumber\\
&  &
\left.-\frac{2 \phi (t)  \phi (1)   \left(\phi^2 (t) - \phi^2 (1)  -r^2\right)  }{\left(\left(\phi (t)- \phi (1 )\right)^2- r^2\right) \left(\left(  \phi (t)+ \phi (1 ) \right)^2-r^2\right)}\right.\nonumber\\
&  &
\times F \left(  i \frac{{\cal{M}}}{1-\ell }+1,  i \frac{{\cal{M}}}{1-\ell };1;\frac{\left( \phi (t)- \phi (1 ) \right) ^2-r^2}{(  \phi (t)+ \phi (1 ) )^2-r^2}\right)
\Bigg]\Bigg| \,dr\,. 
\end{eqnarray*} 
In the notations (\ref{zchange}), it follows  
\begin{eqnarray*}
&  &
\int_0^{\phi (t)- \phi (1 )}   |K_2\left(r,t  ;{\cal{M}};1 \right)| \,dr \\ 
& \leq &
C \left(\left(\phi (t) \right)^2 \right)^{-  i\frac{{\cal{M}}}{1-\ell }} \int_0^{\phi (t)- \phi (1 )}   \Bigg|
\left(\left(z+ 1\right)^2-y^2\right)^{-  i\frac{{\cal{M}}}{1-\ell }}\\
&  &
\times \Bigg[  \frac{ z^2 -z      }{ y^2-\left(z- 1\right)^2 }  
F \left(  i \frac{{\cal{M}}}{1-\ell },  i \frac{{\cal{M}}}{1-\ell };1;\frac{\left(z- 1 \right) ^2-y^2}{( z+ 1 )^2-y^2}\right)\nonumber\\
&  &
\left.-\frac{2 z   \left(1 - z^2    -y^2\right)  }{ \left(\left(z- 1 \right) ^2-y^2\right)
\left(\left(z+ 1 \right) ^2-y^2\right)}\right. 
F \left(  i \frac{{\cal{M}}}{1-\ell }+1,  i \frac{{\cal{M}}}{1-\ell };1;\frac{\left(z- 1 \right) ^2-y^2}{( z+ 1 )^2-y^2}\right)
\Bigg]\Bigg| \,dr \\ 
& \leq &
C \left(\left(\phi (t) \right)^2 \right)^{-  i\frac{{\cal{M}}}{1-\ell }} \int_0^{\phi (t)- \phi (1 )}   \Bigg|
\left(\left(z+ 1\right)^2-y^2\right)^{-  i\frac{{\cal{M}}}{1-\ell }}\frac{  z   }{ \left(z- 1\right)^2- y^2}\\
&  &
\times \Bigg[   -  (z -1)     
F \left(  i \frac{{\cal{M}}}{1-\ell },  i \frac{{\cal{M}}}{1-\ell };1;\frac{\left(z- 1 \right) ^2-y^2}{( z+ 1 )^2-y^2}\right)\nonumber\\
&  &
\left.-\frac{2     \left(1 - z^2    -y^2\right)  }{ \left(z+ 1 \right) ^2-y^2 }\right. 
F \left(  i \frac{{\cal{M}}}{1-\ell }+1,  i \frac{{\cal{M}}}{1-\ell };1;\frac{\left(z- 1 \right) ^2-y^2}{( z+ 1 )^2-y^2}\right)
\Bigg]\Bigg| \,dr \,.
\end{eqnarray*} 
By (\ref{3.20b})
we derive  
\begin{eqnarray*}
&  &
F \left(  i \frac{{\cal{M}}}{1-\ell }+1,  i \frac{{\cal{M}}}{1-\ell };1;\frac{\left( z- 1\right) ^2-y^2}{( z+ 1 )^2-y^2}\right)\\
& = &
- \frac{1}{i \frac{{\cal{M}}}{1-\ell } (\frac{\left( z- 1\right) ^2-y^2}{( z+ 1 )^2-y^2}-1)}\Bigg\{\left(1-i \frac{{\cal{M}}}{1-\ell }\right)   F \left( i \frac{{\cal{M}}}{1-\ell }-1,i \frac{{\cal{M}}}{1-\ell };1;\frac{\left( z- 1\right) ^2-y^2}{( z+ 1 )^2-y^2}\right)\\
&  &
+\left(2 i \frac{{\cal{M}}}{1-\ell }-1\right)  F \left( i \frac{{\cal{M}}}{1-\ell },i \frac{{\cal{M}}}{1-\ell };1;\frac{\left( z- 1\right) ^2-y^2}{( z+ 1 )^2-y^2}\right)\Bigg\} \,.
\end{eqnarray*}
Then
\begin{eqnarray*}
&  &
\int_0^{\phi (t)- \phi (1 )}   |K_2\left(r,t  ;{\cal{M}};1 \right)| \,dr \\ 
& \leq &
C z^{-1 +  2i\frac{{\cal{M}}}{1-\ell }} \int_0^{z- 1}   \Bigg|
\left(\left(z+ 1\right)^2-y^2\right)^{-  i\frac{{\cal{M}}}{1-\ell }} \frac{1}{i 2 {\cal{M}}  \left(y^2-(z-1)^2\right)}\\
&  &
\times \Bigg( i \left(-2 {\cal{M}} \left(y^2+z-1\right)+i (\ell  -1) \left(y^2+z^2-1\right)\right)  F \left(\frac{i {\cal{M}}}{1-\ell  },\frac{i {\cal{M}}}{1-\ell  };1;\frac{(z-1)^2-y^2}{(z+1)^2-y^2}\right)\\
&  &
+ (1-\ell  -i {\cal{M}} ) \left(1-y^2-z^2 \right)  F \left(\frac{i {\cal{M}}}{1-\ell  }-1,\frac{i {\cal{M}}}{1-\ell  };1;\frac{(z-1)^2-y^2}{(z+1)^2-y^2}\right)\Bigg)
\Bigg| \,  d y  \,.
\end{eqnarray*}
For small $\epsilon  $ we have
\begin{eqnarray*}
 F \left(\frac{i {\cal{M}}}{1-\ell },\frac{i {\cal{M}}}{1-\ell  };1;\epsilon \right)
& = & 
1-\frac{{\cal{M}}^2 }{(\ell  -1)^2}\epsilon +\frac{{\cal{M}}^2  (i \ell  +{\cal{M}}-i)^2}{4 (\ell  -1)^4}\epsilon ^2+O\left(\epsilon ^3\right)\,,\\
 F \left(\frac{i {\cal{M}}}{1-\ell }-1 ,\frac{i {\cal{M}}}{1-\ell  };1;\epsilon \right)
& = &
1-\frac{{\cal{M}}   (-i \ell  +{\cal{M}}+i)}{(\ell  -1)^2}\epsilon+\frac{{\cal{M}}^2 
 (-i \ell  +{\cal{M}}+i) (i \ell  +{\cal{M}}-i)}{4 (\ell -1)^4}\epsilon ^2 
+O\left(\epsilon ^3\right)\,.
\end{eqnarray*}
We divide the domain of integration into two zones, according to (\ref{Zone1}) and (\ref{Zone2}) 
and, accordingly,  split the integral into two parts,
\[
\int_0^{z-1}\left|K_2\left(r,t  ;{\cal{M}};1 \right)\right| d r=\int_{(z, y) \in Z_1(\varepsilon, z)}\left|K_2\left(r,t  ;{\cal{M}};1 \right)\right| d y+\int_{(z, y) \in Z_2(\varepsilon, z)}\left|K_2\left(r,t  ;{\cal{M}};1 \right)\right| d y .
\]
In the first zone we have (\ref{1zone}),   
then
\begin{eqnarray*}
&  &
 \frac{1}{i 2 {\cal{M}}  \left(y^2-(z-1)^2\right)}\\
&  &
\times \Bigg( i \left(-2 {\cal{M}} \left(y^2+z-1\right)+i (\ell  -1) \left(y^2+z^2-1\right)\right)  F \left(\frac{i{\cal{M}}}{1-\ell  },\frac{i {\cal{M}}}{1-\ell  };1;\frac{(z-1)^2-y^2}{(z+1)^2-y^2}\right)\\
&  &
+ (1-\ell  -i {\cal{M}} ) \left(1-y^2-z^2 \right)  F \left(\frac{i {\cal{M}}}{1-\ell  }-1,\frac{i {\cal{M}}}{1-\ell  };1;\frac{(z-1)^2-y^2}{(z+1)^2-y^2}\right)\Bigg)
\\
& = & 
-\frac{1}{2} 
+\frac{  \left({\cal{M}}^2 \left(y^2-(z-1)^2\right)+i (\ell  -1) {\cal{M}} \left(y^2+z^2-1\right)+(\ell  -1)^2 \left(y^2+z^2-1\right)\right)}{2 (\ell  -1)^2 \left(y^2-(z-1)^2\right)}\frac{(z-1)^2-y^2}{(z+1)^2-y^2}\\
&  &
-\frac{{\cal{M}} ({\cal{M}}+i (\ell  -1))}{8   (\ell  -1)^4 \left(y^2-(z-1)^2 \right) }\\
&  &
\times \left( {\cal{M}}^2 \left(y^2-(z-1)^2\right)+i (\ell  -1) {\cal{M}} \left(3 y^2+z^2+2 z-3\right)+2 (\ell  -1)^2 \left(y^2+z^2-1\right)\right) \\
&  &
\times \left(\frac{(z-1)^2-y^2}{(z+1)^2-y^2}   \right)^2
 + O\left( \left(\frac{(z-1)^2-y^2}{(z+1)^2-y^2}   \right)^3\right) \,.
\end{eqnarray*}
Hence
\begin{eqnarray*}
&  &
\int_ {Z_1(\epsilon,z)}   |K_2\left(r,t  ;{\cal{M}};1 \right)| \,dr \\ 
& \leq &
C  z  ^{-1+ 2i\frac{{\cal{M}}}{1-\ell }} \int_{Z_1(\epsilon,z)}   \Bigg|
\left(\left(z+ 1\right)^2-y^2\right)^{-  i\frac{{\cal{M}}}{1-\ell }} \nonumber \\
&  &
\times
\Bigg[-\frac{1}{2} \nonumber 
+\frac{  \left({\cal{M}}^2 \left(y^2-(z-1)^2\right)+i (\ell  -1){\cal{M}} \left(y^2+z^2-1\right)+(\ell  -1)^2 \left(y^2+z^2-1\right)\right)}{2 (\ell -1)^2 \left(y^2-(z-1)^2\right)}\frac{(z-1)^2-y^2}{(z+1)^2-y^2} \nonumber \\
&  &
-\frac{{\cal{M}} ({\cal{M}}+i (\ell  -1))}{8   (\ell  -1)^4 \left(y^2-(z-1)^2\right)  } \nonumber \\
&  &
\times \left( {\cal{M}}^2 \left(y^2-(z-1)^2\right)+i (\ell  -1) {\cal{M}} \left(3 y^2+z^2+2 z-3\right)+2 (\ell  -1)^2 \left(y^2+z^2-1\right)\right) \nonumber\\
&  &
\times \left(\frac{(z-1)^2-y^2}{(z+1)^2-y^2}   \right)^2 \nonumber 
 + O\left( \left(\frac{(z-1)^2-y^2}{(z+1)^2-y^2}   \right)^3\right)\Bigg]\Bigg| dy \,. 
 \end{eqnarray*}
 We split integral into three parts and consider the first one.  Since $i\frac{{\cal{M}}}{1-\ell } <1$, we obtain 
 \begin{eqnarray*}
A_0
&:=  & 
 z  ^{-1+  2i\frac{{\cal{M}}}{1-\ell }}\int_ {Z_1(\epsilon,z)}   
\left(\left(z+ 1\right)^2-y^2\right)^{-  i\frac{{\cal{M}}}{1-\ell }}  dy\\
&\lesssim & 
 z  ^{-1+  2i\frac{{\cal{M}}}{1-\ell }}\int_0^{z-1}  
\left(\left(z+ 1\right)^2-y^2\right)^{-  i\frac{{\cal{M}}}{1-\ell }}  dy\\
& = &
 z  ^{-1+  2i\frac{{\cal{M}}}{1-\ell }}(z-1) (z+1)^{ - 2 i\frac{{\cal{M}}}{1-\ell } }  F \left(\frac{1}{2},  i\frac{{\cal{M}}}{1-\ell } ;\frac{3}{2};\frac{(z-1)^2}{(z+1)^2}\right)\\\\
& \lesssim &  
1  +  \left(   \frac{4z }{ (z+ 1  )^2} \right)^{ 1 +  i \frac{{\cal{M}}}{\ell-1}}.
 \end{eqnarray*}
Hence $A_0 \lesssim 1$ for all $z \in [1,\infty)$.  Similarly,  
\begin{eqnarray*} 
A_1 
& := &
  z  ^{-1+ 2i\frac{{\cal{M}}}{1-\ell }} \int_ {Z_1(\epsilon,z)}  \frac{(z-1)^2-y^2}{(z+1)^2-y^2}
 \Bigg|
\left(\left(z+ 1\right)^2-y^2\right)^{-  i\frac{{\cal{M}}}{1-\ell }}\\
&  &
\times\Bigg[\frac{  \left({\cal{M}}^2 \left(y^2-(z-1)^2\right)+i (\ell -1) {\cal{M}} \left(y^2+z^2-1\right)+(\ell  -1)^2 \left(y^2+z^2-1\right)\right)}{2 (\ell  -1)^2 \left(y^2-(z-1)^2\right)}\Bigg]\Bigg|dy\\
& \lesssim  &
  z  ^{-1+ 2i\frac{{\cal{M}}}{1-\ell }} \int_0^{z-1}   \Bigg|
\left(\left(z+ 1\right)^2-y^2\right)^{-  i\frac{{\cal{M}}}{1-\ell }}\frac{ z^2}{\left(y^2-(z-1)^2\right)}\frac{(z-1)^2-y^2}{(z+1)^2-y^2}
 \Bigg|dy\\
& \lesssim  &
  z  ^{ 1+ 2i\frac{{\cal{M}}}{1-\ell }} \int_0^{z-1}  
\left(\left(z+ 1\right)^2-y^2\right)^{-1-  i\frac{{\cal{M}}}{1-\ell }}
 dy\\
& \lesssim  &
    F \left(\frac{1}{2},1+ i\frac{{\cal{M}}}{1-\ell } ;\frac{3}{2};\frac{(z-1)^2}{(z+1)^2}\right) \,.
 \end{eqnarray*}
  According to 
Lemma~\ref{L4.1}  
 we have 
 \begin{eqnarray*} 
 A_1  
& \lesssim  &
\cases{  z  ^{ \frac{1}{  2}-\frac{1}{  2(\ell - 1 ) } \sqrt{(1-3 \ell )^2-4 m_c^2} )}  \quad  
 \mbox{if }\quad 4m_c^2>  4\ell(2\ell-1)\,,   \cr
1 \quad   \mbox{if }\quad 4m_c^2\leq 4\ell(2\ell-1) \,. }
 \end{eqnarray*} 
Similarly,  
\begin{eqnarray*}
A_2 
& := &
  z  ^{-1+ 2i\frac{{\cal{M}}}{1-\ell }} \int_ {Z_1(\epsilon,z)}    \Bigg|
\left(\left(z+ 1\right)^2-y^2\right)^{-  i\frac{{\cal{M}}}{1-\ell }}\\
&  &
\times\Bigg[
-\frac{{\cal{M}} ({\cal{M}}+i (\ell -1))}{8   (\ell -1)^4 \left(y^2-(z-1)^2\right) }\\
&  &
\times \Big({\cal{M}}^2 \left(y^2-(z-1)^2\right)+i (\ell -1) {\cal{M}} \left(3 y^2+z^2+2 z-3\right)+2 (\ell -1)^2 \left(y^2+z^2-1\right) \Big) \\
&  &
\times \left( \frac{(z-1)^2-y^2}{(z+1)^2-y^2}\right)^2\Bigg]\Bigg|dy\\
& \lesssim  &
  z  ^{-1+ 2i\frac{{\cal{M}}}{1-\ell }} \int_0^{z-1}   \Bigg|
\left(\left(z+ 1\right)^2-y^2\right)^{-  i\frac{{\cal{M}}}{1-\ell }}\Bigg[
\frac{1}{ \left(y^2-(z-1)^2\right) } z^2 \left( \frac{(z-1)^2-y^2}{(z+1)^2-y^2}\right)^2\Bigg]\Bigg|dy \\
& \lesssim  &
 z  ^{1+ 2i\frac{{\cal{M}}}{1-\ell }}(z-1) (z+1)^{2(-1-  i\frac{{\cal{M}}}{1-\ell })}  
F \left(\frac{1}{2}, 1+ i\frac{{\cal{M}}}{1-\ell };\frac{3}{2};\frac{(z-1)^2}{(z+1)^2}\right)\\
& \lesssim  &
\cases{ z  ^{ \frac{1}{  2}-\frac{1}{  2(\ell - 1 ) } \sqrt{(1-3 \ell )^2-4 m_c^2} )}  \quad   \mbox{if }\quad 4m_c^2>  4\ell(2\ell-1)\,,   \cr
1 \quad   \mbox{if }\quad 4m_c^2\leq 4\ell(2\ell-1) \,. }
 \end{eqnarray*}
 The last term is defined by
\begin{eqnarray*}
A_3 
& := &
z  ^{-1+ 2i\frac{{\cal{M}}}{1-\ell }} \int_{Z_1(\epsilon,z)}   \Bigg|
\left(\left(z+ 1\right)^2-y^2\right)^{-  i\frac{{\cal{M}}}{1-\ell }} 
\Bigg[ O\left( \left(\frac{(z-1)^2-y^2}{(z+1)^2-y^2}   \right)^3\right)\Bigg] \Bigg| dy 
 \end{eqnarray*}
and estimated as follows:
\begin{eqnarray*}
A_3 
& \lesssim &
z  ^{-1+ 2i\frac{{\cal{M}}}{1-\ell }} \int_0^{z-1}   
\left(\left(z+ 1\right)^2-y^2\right)^{-  i\frac{{\cal{M}}}{1-\ell }}\frac{(z-1)^2-y^2}{(z+1)^2-y^2}  
 dy \\
& \lesssim & 
\cases{ z  ^{ \frac{1}{  2}-\frac{1}{  2(\ell - 1 ) } \sqrt{(1-3 \ell )^2-4 m_c^2} )}  \quad   \mbox{if }\quad 4m_c^2>  4\ell(2\ell-1)\,,   \cr
1 \quad   \mbox{if }\quad 4m_c^2\leq 4\ell(2\ell-1) \,.} 
 \end{eqnarray*}

In $Z_2(\varepsilon,z) $ 
we have (\ref{2zone})  
and 
\begin{eqnarray*}
\left|  F \left(\frac{i {\cal{M}}}{1-\ell },\frac{i {\cal{M}}}{1-\ell  };1;\frac{(z-1)^2-y^2}{(z+1)^2-y^2}\right) \right| 
& \leq & C\quad \mbox{for all}\quad z \in[1,\infty),\quad y \in [0,z-1]\,, \\
\left|  F \left(\frac{i {\cal{M}}}{1-\ell  }-1,\frac{i {\cal{M}}}{1-\ell  };1;\frac{(z-1)^2-y^2}{(z+1)^2-y^2}\right) \right| 
& \leq & C \quad \mbox{for all}\quad z \in[1,\infty),\quad y \in [0,z-1]\,,
 \end{eqnarray*}
then
\begin{eqnarray*}
&  &
\int_ {Z_2(\epsilon,z)}   |K_2\left(r,t  ;{\cal{M}};1 \right)| \,dr \\ 
& \lesssim  &
 z^{-1 +  2i\frac{{\cal{M}}}{1-\ell }} \int_ {Z_2(\epsilon,z)}  \Bigg|
\left(\left(z+ 1\right)^2-y^2\right)^{-  i\frac{{\cal{M}}}{1-\ell }}\Bigg[\frac{1}{i 2 {\cal{M}}  \left(y^2-(z-1)^2\right)}\\
&  &
\times\Bigg(i \left(-2 {\cal{M}} \left(y^2+z-1\right)+i (\ell  -1) \left(y^2+z^2-1\right)\right)  F \left(\frac{i {\cal{M}}}{1-\ell },\frac{i {\cal{M}}}{1-\ell  };1;\frac{(z-1)^2-y^2}{(z+1)^2-y^2}\right)\\
&  &
+ (1-\ell  -i {\cal{M}} ) \left(1-y^2-z^2 \right)  F \left(\frac{i {\cal{M}}}{1-\ell  }-1,\frac{i {\cal{M}}}{1-\ell  };1;\frac{(z-1)^2-y^2}{(z+1)^2-y^2}\right)\Bigg)
\Bigg]\Bigg| \,  d y  \\ 
& \lesssim  &
 z^{-1 +  2i\frac{{\cal{M}}}{1-\ell }} \int_0^{z- 1}   \Bigg|
\left(\left(z+ 1\right)^2-y^2\right)^{-  i\frac{{\cal{M}}}{1-\ell }}\Bigg[\frac{1}{i 2 {\cal{M}}  \left(y^2-(z-1)^2\right)}z^2 
\Bigg]\Bigg| \,  d y  \\ 
& \lesssim  &
 \varepsilon^{-1}z^{ 1 +  2i\frac{{\cal{M}}}{1-\ell }} \int_0^{z- 1}   
\left(\left(z+ 1\right)^2-y^2\right)^{-1-  i\frac{{\cal{M}}}{1-\ell }}   \,  d y  \\ 
& \lesssim  &
\cases{ z  ^{ \frac{1}{  2}-\frac{1}{  2(\ell - 1 ) } \sqrt{(1-3 \ell )^2-4 m_c^2} )}  \quad   \mbox{if }\quad 4m_c^2>  4\ell(2\ell-1) ,  \cr
1 \quad   \mbox{if }\quad 4m_c^2\leq 4\ell(2\ell-1) , }
\end{eqnarray*}
 for all $ z \in[1,\infty)$. The lemma is proved. \qed

\subsection{Problem  without  potential and nonlinear term. Proof of Theorem~\ref{T_small_mass}} 
\label{SS4.3}

\begin{proposition}
\label{P4.1}
Let $\ell>1$ and $\gamma \in \R$. Then the solution $\psi_{ID}(x,t )$ of the linear equation (\ref{3.23}) that takes initial values (\ref{13.12}) satisfies the estimate
\begin{eqnarray*}
t^\gamma \| \psi_{ID}(x,t )\|_{H_{(s)}}
& \leq &
Ct^{\gamma +k_+}( \| \psi _0 \|_{H_{(s)}}+ \| \psi _1 \|_{H_{(s)}}) 
+
Ct^{\gamma  - \ell  } 
\| \psi _0  \|_{H_{(s)}}\\
&  &
+Ct^{\gamma +k_+}
\|\psi _0 \|_{H_{(s)}}\cases{  t  ^{ \frac{  \ell - 1 }{  2}-\frac{1}{  2 } \sqrt{(1-3 \ell )^2-4 m_c^2} )}  \quad   \mbox{if }\quad 4m_c^2>  4\ell(2\ell-1)  , \cr
1 \quad   \mbox{if }\quad 4m_c^2\leq 4\ell(2\ell-1)  }  
\end{eqnarray*}
for all $t \in [1,\infty)$.
\end{proposition}
\msk

\ndt
{\bf Proof.} 
The solution, after the partial Liouville transformation,  can be written as follows:
\begin{eqnarray*}
 \psi_{ID}(x,t )
& = &
t^{k_+}\int_0^{\phi (t)- \phi (1 )}  K_1 \left(r,t; {\cal{M}} ;1 \right) v_{u _1}(x,r)\,dr 
+
t^{k_+}\left(\frac{\phi (t)}{ \phi (1)} \right)^{-   i\frac{{\cal{M}}}{1-\ell} } 
v_{u _0}\left(x,\phi (t)- \phi (1) \right)\\
&  &
+t^{k_+}\frac{1}{\phi (1)}\int_0^{\phi (t)- \phi (1)}  
K_2\left(r,t  ;{\cal{M}};1 \right)
v_{u _0}(x, r)\,dr  \,.
\end{eqnarray*}
Then, by (\ref{Blplq}), we obtain 
\begin{eqnarray*}
\| \psi_{ID}(x,t )\|_{H_{(s)}}
& \lesssim &
t^{k_+}\|u _1  \|_{H_{(s)}}\int_0^{\phi (t)- \phi (1 )}   | K_1 \left(r,t; {\cal{M}} ;1 \right)|  \,dr
+
t^{k_+}\left(\frac{\phi (t)}{ \phi (1)} \right)^{-   i\frac{{\cal{M}}}{1-\ell} } 
\| u _0  \|_{H_{(s)}}\\
&  &
+t^{k_+}
\|u _0 \|_{H_{(s)}}\int_0^{\phi (t)- \phi (1)}  
|K_2\left(r,t  ;{\cal{M}};1 \right)|\,dr \\
&  \lesssim  &
t^{k_+}\|u _1  \|_{H_{(s)}} 
+
t^{ - \ell  } 
\| u _0  \|_{H_{(s)}}\\
&  &
+t^{k_+}
\|u _0 \|_{H_{(s)}}
\cases{  t  ^{ \frac{ (\ell - 1 )}{  2}-\frac{1}{  2 } \sqrt{(1-3 \ell )^2-4 m_c^2} )}  \quad   \mbox{if }\quad 4m_c^2>  4\ell(2\ell-1)   ,\cr
1 \quad   \mbox{if }\quad 4m_c^2\leq 4\ell(2\ell-1)  \,.} 
\end{eqnarray*}
Here 
$u _0(x, t)=\psi_0(x, t),\quad  u _1(x, t)=\psi_1(x, t)- k_+\psi_0(x, t)$. 
Hence, the proposition is proved. \qed

\subsection{Problem with vanishing initial conditions and without  potential. Estimate of nonlinear term. Proof of Theorem~\ref{T_small_mass}}

In this subsection, we obtain a 
decay estimates of the solution of the linear equation 
\begin{eqnarray*}
    \frac{\partial^2 \psi }{\partial t^2}
+ \frac{3 \ell}{t}    \frac{\partial \psi }{\partial t}
-  t^{-2\ell}{\mathcal{A}}(x,\partial_x)\psi +  t^{ - 2 }  \frac{m^2c^4}{h^2} \psi = \Psi\,.
\end{eqnarray*}
For $    k_ +   = - \frac{1}{2}( 3 \ell - 1- \sqrt{(3 \ell - 1 )^2 -4 m_c^2} )$, we define operator $G$   by 
\begin{eqnarray*}
G[F]
& :=  &
  t^{k_+}\int_1 ^{t  } b^{-k_+} \,d b\int_0^{ \phi (t)- \phi (b ) }  
 E(r,t;b ; {\cal{M}}) {\cal E}{\cal E}[F]  \left(x,r;b  \right) \,dr\,.
\end{eqnarray*}
Hence, the operator $G$  gives a solution to the problem
\begin{eqnarray*}
\cases{   \dsp   \frac{\partial^2 \psi }{\partial t^2}
+ \frac{3 \ell}{t}    \frac{\partial \psi }{\partial t}
-  t^{-2\ell}{\mathcal{A}}_{3/2}(x,\partial_x)\psi +  t^{ - 2 }   m_c^2  \psi =F \,,\cr
\dsp \psi (x,1)=0,\quad \psi_t (x,1)=0.}
\end{eqnarray*}
 
\begin{theorem}
\label{T4.7}
Assume that   $\ell >1 $.  
Then for the operator $G$   the following estimate holds 
\begin{eqnarray*} 
t^{\gamma}\| G[ \Psi(\Phi)]\|_{H_{(s)}} 
& \leq &
C_{NT} \left(\max_{b \in [1,t]}  \tau^{\gamma} \| \Phi   \left(x, \tau  \right)\| _{H_{(s)}}\right)^{1+\alpha} \\
& &
\times \cases{ t^{ 2-\gamma \alpha } \quad \mbox{if} \quad  \gamma< \frac{3\ell+3-\sqrt{(1-3 \ell )^2-4 m_c^2}}{2(1+\alpha)} ,  \cr 
  t^{\gamma+\frac{1}{2}\left( 1-3\ell +\sqrt{(1-3 \ell )^2-4 m_c^2}\right)} 
\quad \mbox{if} \quad \gamma>\frac{3\ell+3-\sqrt{(1-3 \ell )^2-4 m_c^2}}{2(1+\alpha)},\cr
t^{\gamma+\frac{1}{2}\left( 1-3\ell +\sqrt{(1-3 \ell )^2-4 m_c^2}\right)} \ln (t)\quad 
\mbox{if} \quad \gamma=\frac{3\ell+3-\sqrt{(1-3 \ell )^2-4 m_c^2}}{2(1+\alpha)} }    
\end{eqnarray*}
for all $t \in [1,\infty)$.
\end{theorem}
\msk

\ndt
{\bf Proof.} Indeed,  from (\ref{G}), (\ref{Blplq}),  and Lemma~\ref{L4.2E} we obtain
\begin{eqnarray*} 
\| G[F]\|_{ H_{(s)}} 
 & \lesssim &
 t^{k_+} \int_1 ^{t  } b^{-k_+} \| F    (x,b)\|_{ H_{(s)}}\,d b\int_0^{ \phi (t)- \phi (b ) }  
 |E(r,t;b ; {\cal{M}})|  \,dr\\ 
 &  \lesssim  &
 t^{k_+} \int_1 ^{t  } b^{-k_+} \| F    (x,b)\|_{ H_{(s)}}
 | \phi (b ) |^{-\frac{1}{\ell-1}-\frac{\sqrt{(1-3 \ell )^2-4 m_c^2}}{\ell-1}}  
   | \phi (b)+\phi (t) |^{ \frac{\sqrt{(1-3 \ell )^2-4 m_c^2}}{ \ell-1 }}\,d b\\  
 &  \lesssim &
 t^{k_+} \int_1 ^{t  } b^{-k_+ +1}  \| F    (x,b)\|_{ H_{(s)}}\,d b\,.
\end{eqnarray*}
Further, \begin{eqnarray*}  
t^\gamma\| G[ \Psi(\Phi)]\|_{H_{(s)}} 
& \lesssim &
 t^{\gamma+k_+} \int_1 ^{t  } b^{-k_+ +1}    \left(\| \Phi   \left(x, b  \right)\| _{H_{(s)}}\right)^{1+\alpha}\,d b\\ 
 & \lesssim &
\left(\max_{b \in [1,t]}  b^{\gamma} \| \Phi   \left(x, b  \right)\| _{H_{(s)}}\right)^{1+\alpha} I (t)\,,
\end{eqnarray*}
where with $2-k_+  -\gamma(1+\alpha)\not= 0 $ we obtain 
\begin{eqnarray*} 
I(t)
& := &
t^{\gamma+k_+} \int_1 ^{t  } b^{-k_+ +1-\gamma(1+\alpha)}  \,d b\\
& \lesssim &
\cases{ t^{2 -\gamma \alpha } \quad \mbox{if} \quad  \gamma< \frac{3\ell+3-\sqrt{(1-3 \ell )^2-4 m_c^2}}{2(1+\alpha)}  \,, \cr 
  t^{\gamma+\frac{1}{2}\left( 1-3\ell +\sqrt{(1-3 \ell )^2-4 m_c^2}\right)} \quad \mbox{if} \quad \gamma>\frac{3\ell+3-\sqrt{(1-3 \ell )^2-4 m_c^2}}{2(1+\alpha)} \,.}    
\end{eqnarray*}
For $2-k_+  -\gamma(1+\alpha) = 0 $ the estimate 
$
I(t)
  \lesssim  
t^{\gamma+k_+} \ln (t)    
$ is evident. 
The theorem is proved. \qed \\

\subsection{Estimate of equation with potential and with vanishing initial data.  Proof of Theorem~\ref{T_small_mass}}

\begin{theorem}
\label{T4.8}
Assume that with $\delta  \geq 2 $ the potential  $V$ satisfies the condition (\ref{V}).  
Then
\begin{eqnarray*}
t^{\gamma} \| G[V\Phi]\|_{ H_{(s)}}
& \leq & 
C \varepsilon_P   \sup_{\tau \in [1,t]} \tau^{\gamma} \| \Phi    (x,\tau)\|_{ H_{(s)}}  \\
& &
\times
\cases{ t^{2-\delta}\quad \mbox{if} \quad  \gamma< \frac{3\ell+3-\sqrt{(1-3 \ell )^2-4 m_c^2}}{2 }-\delta,    \cr 
  t^{\gamma+\frac{1}{2}\left( 1-3\ell +\sqrt{(1-3 \ell )^2-4 m_c^2}\right)} \quad \mbox{if} \quad \gamma>\frac{3\ell+3-\sqrt{(1-3 \ell )^2-4 m_c^2}}{2 }-\delta , \cr
  t^{\gamma+\frac{1}{2}\left( 1-3\ell +\sqrt{(1-3 \ell )^2-4 m_c^2}\right)} \ln (t) \quad \mbox{if} \quad \gamma=\frac{3\ell+3-\sqrt{(1-3 \ell )^2-4 m_c^2}}{2 } -\delta }  
\end{eqnarray*}
 for all $ t \in[1,\infty)$.  
\end{theorem}  
\medskip

\ndt
{\bf Proof.} Indeed, from (\ref{G}), (\ref{Blplq}), and Lemma~\ref{L4.2E}  we obtain
\begin{eqnarray*} 
t^\gamma \| G[V\Phi]\|_{ H_{(s)}} 
 & \leq &
 t^{\gamma +k_+} \int_1 ^{t  } b^{-k_+} \,d b\int_0^{ \phi (t)- \phi (b ) }  
 |E(r,t;b ; {\cal{M}})| \|{\cal E E} [V\Phi]   (x,r;b)\|_{ H_{(s)}} \,dr\\
 & \leq &
 t^{\gamma +k_+} \int_1 ^{t  } b^{-k_+} \| V\Phi    (x,b)\|_{ H_{(s)}}\,d b\int_0^{ \phi (t)- \phi (b ) }  
 |E(r,t;b ; {\cal{M}})|  \,dr\\
& \leq &
 C \varepsilon_P  t^{\gamma + k_+} \int_1 ^{t  } b^{-k_+-\delta} \| \Phi    (x,b)\|_{ H_{(s)}}  | \phi (b )^{-\frac{1}{\ell-1}-\frac{\sqrt{(1-3 \ell )^2-4 m_c^2}}{\ell-1}}  
   ( \phi (b)+\phi (t) )^{ \frac{\sqrt{(1-3 \ell )^2-4 m_c^2}}{ \ell-1 }} |\,d b \\
& \leq &
 C \varepsilon_P  t^{\gamma + k_+} \int_1 ^{t  } b^{-k_++1-\delta} \| \Phi    (x,b)\|_{ H_{(s)}}  \,d b \\
& \leq &
 C \varepsilon_P  \left( \sup_{\tau \in [1,t]} \tau^{\gamma} \| \Phi    (x,\tau)\|_{ H_{(s)}} \right)J(t)  \,,
\end{eqnarray*} 
where
$
J (t)  := 
t^{\gamma + k_+} \int_1 ^{t  } b^{-k_++1-\delta -\gamma}    \,d b$. 
With $2-k_+ -\delta -\gamma \not= 0 $ we obtain 
\begin{eqnarray*} 
J(t)
& \lesssim &
\cases{ t^{2-\delta}\quad \mbox{if} \quad  \gamma< \frac{3\ell+3-\sqrt{(1-3 \ell )^2-4 m_c^2}}{2 }-\delta ,  \cr 
  t^{\gamma+\frac{1}{2}\left( 1-3\ell +\sqrt{(1-3 \ell )^2-4 m_c^2}\right)} 
\quad \mbox{if} \quad \gamma>\frac{3\ell+3-\sqrt{(1-3 \ell )^2-4 m_c^2}}{2 }-\delta  . }    
\end{eqnarray*}
If $2-k_+  -\gamma -\delta   = 0 $, then 
$ J(t)    \lesssim  t^{\gamma + k_+} \ln (t)$.
  The theorem is proved.  
\qed

\subsection{Existence of solution to semilinear equation.  Proof of Theorem~\ref{T_small_mass}}
\label{SS4.6}

We solve the Cauchy problem    (\ref{KGE_KF})\&(\ref{13.12})
    through the integral equation.
We  appeal to the operator (\ref{G}),  
$
G:={\mathcal K}\circ {\mathcal EE} 
$, 
where ${\mathcal EE}$ stands for the evolution (wave) equation in the exterior of BH in the universe without expansion, while ${\mathcal K}$ is introduced  by (\ref{9.9}) with $k_+ $, ${\cal{M}_+} $. The kernel $ E(r,t; 0,b;M) $ is given by  (\ref{Edef}). For the function $f(x,t) $ we define
$
v(x,t;b):= {\mathcal EE} [f](x,t;b)$, 
where the function 
$v(x,t;b)$   
is a solution to the Cauchy problem (\ref{2.15}).

We intend to apply Banach's fixed-point theorem.  
We use the Lipschitz condition (${\mathcal L}$) to estimate
nonlinear term.
We start with the integral equation (\ref{5.1}),   
where the function $ \psi _{ID} (x,t) \in C([0,\infty);H_{(s)})$ is given.

Consider the mapping $S [\Phi]$ (\ref{14.5}), 
where the function $
\psi_{ID}  
$ 
is generated by initial data, that is, by  (\ref{13.12}). 
It is easily seen that   supp\,$S [\Phi] \subset$ supp\,$ \Phi $ if   supp\,$ \Phi  \subset$ supp\,$ \psi_{ID} $.   We claim that if $ \Phi \in X({R,H_{(s)},\gamma})$ with   $\gamma \in [0,H] $ and if  supp~$\Phi \subseteq \{(x,t) \in \R^3\times[0,\infty)\,|\,|x|> R_{ID}- ( \phi (t)-\phi(1)) \,\} $, then 
$S[\Phi] \in   X({R,H_{(s)},\gamma})$.  Moreover, $S$ 
is a contraction, provided that   $\varepsilon  $, $\varepsilon_P  $, and $R$ are sufficiently small.

First, we note that due to Proposition~\ref{P4.1}, we have
\begin{eqnarray*}
t^\gamma \| \psi_{ID}(x,t )\|_{H_{(s)}}
& \leq &
C_{ID}t^{\gamma +k_+}( \| \psi _0 \|_{H_{(s)}}+ \| \psi _1 \|_{H_{(s)}}) 
+
C_{ID}t^{\gamma  - \ell  } 
\| \psi _0  \|_{H_{(s)}}\\
&  &
+C_{ID}t^{\gamma +k_+}
\|\psi _0 \|_{H_{(s)}}\cases{ t  ^{ \frac{  \ell - 1 }{  2}-\frac{1}{  2 } \sqrt{(1-3 \ell )^2-4 m_c^2} )}  \quad   
\mbox{if }\quad 4m_c^2>  4\ell(2\ell-1)  , \cr
1 \quad   \mbox{if }\quad 4m_c^2\leq 4\ell(2\ell-1) \,. } 
\end{eqnarray*}
According to Theorem~\ref{T4.8}
\begin{eqnarray*}
t^{\gamma} \| G[V\Phi]\|_{ H_{(s)}}
& \leq & 
C \varepsilon_P   \sup_{\tau \in [1,t]} \tau^{\gamma} \| \Phi    (x,\tau)\|_{ H_{(s)}}  \\
& &
\times
\cases{ t^{2-\delta}\quad \mbox{if} \quad  \gamma< \frac{3\ell+3-\sqrt{(1-3 \ell )^2-4 m_c^2}}{2 }-\delta,    \cr 
  t^{\gamma+\frac{1}{2}\left( 1-3\ell +\sqrt{(1-3 \ell )^2-4 m_c^2}\right)} \quad \mbox{if} \quad \gamma>\frac{3\ell+3-\sqrt{(1-3 \ell )^2-4 m_c^2}}{2 }-\delta , \cr
  t^{\gamma+\frac{1}{2}\left( 1-3\ell +\sqrt{(1-3 \ell )^2-4 m_c^2}\right)}\ln (t) \quad \mbox{if} \quad \gamma=\frac{3\ell+3-\sqrt{(1-3 \ell )^2-4 m_c^2}}{2 } -\delta ,}  
\end{eqnarray*} 
while due to Theorem~\ref{T4.7} \begin{eqnarray*} 
t^{\gamma}\| G[ \Psi(\Phi)]\|_{H_{(s)}} 
& \leq &
C_{NT} \left(\max_{b \in [1,t]}  \tau^{\gamma} \| \Phi   \left(x, \tau  \right)\| _{H_{(s)}}\right)^{1+\alpha} \\
& &
\times \cases{ t^{2 -\gamma \alpha } \quad \mbox{if} \quad  \gamma< \frac{3\ell+3-\sqrt{(1-3 \ell )^2-4 m_c^2}}{2(1+\alpha)}  , \cr 
  t^{\gamma+\frac{1}{2}\left( 1-3\ell +\sqrt{(1-3 \ell )^2-4 m_c^2}\right)} \quad \mbox{if} \quad \gamma>\frac{3\ell+3-\sqrt{(1-3 \ell )^2-4 m_c^2}}{2(1+\alpha)}, \cr
t^{\gamma+\frac{1}{2}\left( 1-3\ell +\sqrt{(1-3 \ell )^2-4 m_c^2}\right)} \ln (t)\quad \mbox{if} \quad \gamma=\frac{3\ell+3-\sqrt{(1-3 \ell )^2-4 m_c^2}}{2(1+\alpha)} .}    
\end{eqnarray*}

  The condition (\ref{V}) and the estimate for the potential term with $\delta =2$ lead  to the  choice \\
 $\gamma< \frac{3\ell-1-\sqrt{(1-3 \ell )^2-4 m_c^2}}{2}  $. Hence, the time-dependent factor  of the potential term
is bounded.  
Further, the time-dependent factor  of the nonlinear term decays for $\gamma< \frac{3\ell+3-\sqrt{(1-3 \ell )^2-4 m_c^2}}{2(1+\alpha)}$ and $2-\gamma \alpha <0$, which imply $\alpha >\frac{4}{3\ell-1-\sqrt{(1-3 \ell )^2-4 m_c^2}} $. The last inequality  is equivalent to 
\begin{equation}
\label{4.32}
\frac{3\ell-1-\sqrt{(1-3 \ell )^2-4 m_c^2}}{2} > \frac{3\ell+3-\sqrt{(1-3 \ell )^2-4 m_c^2}}{2(1+\alpha)}\,.
\end{equation} 
It is easy to see that if (\ref{4.32}) is violated, then none of the estimates for the nonlinear term gives a bounded factor. 
The condition (\ref{4.32})  allows us to   choose  
$
\gamma \in \left( \frac{2}{\alpha}, \frac{3\ell+3-\sqrt{(1-3 \ell )^2-4 m_c^2}}{2(1+\alpha)} \right) 
$ in the first case of the estimate of nonlinear factor or   $
\gamma \in \left[  
\frac{3\ell+3-\sqrt{(1-3 \ell )^2-4 m_c^2}}{2(1+\alpha)} ,  \frac{3\ell-1-\sqrt{(1-3 \ell )^2-4 m_c^2}}{2} \right)$ due to the remaining two cases. 

If $\delta >2$,  then the   time-dependent factor  of the  potential term  decays. The estimates for the nonlinear term due to the condition on $\gamma$ in Theorem~\ref{T_small_mass} lead  to a bounded factor.

 Thus, the last inequalities prove  that the operator $S$ maps $X({R,s,\gamma})$ into itself if 
$
\| \psi _0 \|_{H_{(s)}}+ \| \psi _1 \|_{H_{(s)}} \leq \varepsilon
$, while 
$\varepsilon_P  $, $\varepsilon  $, and $R$ are sufficiently small, namely, provided that
$
\frac{1}{ 1- \varepsilon_P C_{P} } C_{ID} 
\varepsilon  + \frac{1}{ 1- \varepsilon_P C_{P}  }C_{NT} R^{\alpha +1} < R  
$ and $\gamma < \frac{3\ell-1-\sqrt{(1-3 \ell )^2-4 m_c^2}}{2} $.
To prove that $S$ is a contraction mapping, we derive the contraction property from
\[
\sup_{t \in [0,\infty) }t^{ \gamma } \|S[\Phi_1 ](\cdot,t) -  S[\Phi_2  ](\cdot,t) \|_{H_{(s)}({\mathbb R}^n) }
 \leq  
CR ^{\alpha } d(\Phi_1,\Phi_2 )\,.  
\]
The theorem is proved. \qed

\section*{Acknowledgments}

The authors are indebted to the referees for  the  
valuable remarks  and suggestions, which improved the text.


\begin{thebibliography}{99}




\bibitem{Aslan}
H. S. Aslan, M. R. Ebert and M. Reissig, {\it Scale--invariant  semilinear damped wave models with mass term and integrable in time speed of propagation},   Differential and Integral Equations, 36 (2023), 453--490.

 

\bibitem{B-E}
H. Bateman and  A. Erdelyi, ``Higher  Transcendental  Functions",  Vol. 1, 2, McGraw-Hill, 1954. 



\bibitem{Gibbons}
 S. Chen, G. W. Gibbons, Y. Li and  Y. Yang, {\it   Friedmann's equations in all dimensions and Chebyshev's theorem}, J. Cosmol. Astropart. Phys., 12 (2014), 035.


\bibitem{Costa}  
J. L. Costa, A. T. Franzen and J. Oliver, {\it Semilinear wave equations on accelerated expanding FLRW spacetimes}, { Ann. Henri Poincar\'e},  24 (2023), 3185-3207.


 \bibitem{Farrah}
D. Farrah et al., {\it Observational evidence for cosmological coupling of black holes and its implications for an astrophysical source of dark energy}, {The Astrophysical Journal Letters}, 944 (2023), L31.
 
 

\bibitem{NA2015}
A. Galstian and  K. Yagdjian, {\it Global solutions for semilinear Klein-Gordon equations in FLRW spacetimes},  Nonlinear Anal., 113 (2015), 339-356. 


 \bibitem{Nakamura2014}
M. Nakamura, {\it The Cauchy problem for semi-linear Klein-Gordon equations in de Sitter spacetime}, J. Math. Anal. Appl., 410 (2014), 445-454.
  
 \bibitem{NatarioAHP2023}
J. Nat\'ario and F. Rossetti, {\it Explicit formulas and decay rates for the solution of the wave equation in cosmological spacetimes}, J. Math. Phys., 64 (2023), 032504.
   
   
\bibitem{NatarioAHP2022}
J. Nat\'ario and A. Sasane, {\it Decay of solutions to the Klein-Gordon equation on some expanding cosmological spacetimes},  Ann. Henri Poincar\'e, 23 (2022), 2345--2389.

 \bibitem{Palmieri}
A. Palmieri, {\it  Integral representation formulae for the solution of a wave equation with time-dependent damping and mass in the scale-invariant case}, 
Math. Methods Appl. Sci., {44} (2021), 13008-13039.


\bibitem{Perlick}
 V. Perlick, O. Y. Tsupko and G. S. Bisnovatyi-Kogan, {\it Black hole shadow in an expanding universe with a cosmological constant}, Phys. Rev. D, 97 (2018), 104062.

 
\bibitem{Taylor}   
M. Taylor, ``Partial Differential Equations III. Nonlinear Equations", 2nd edition, Applied Mathematical Sciences  117, Springer, New York, 2011.

\bibitem{Tsutaya-Wakasugi} 
K.Tsutaya,  Y. Wakasugi,  {\it Blow-up of solutions of semilinear wave equations in accelerated expanding Friedmann-Lema\^itre-Robertson-Walker spacetime}, Rev. Math. Phys. 34 (2022), no. 2, Paper No. 2250003, 16 pp. 
 
\bibitem{YagBook}
K. Yagdjian,   
`` The Cauchy Problem for Hyperbolic Operators.  
Multiple Characteristics. Micro-Local Approach",  Mathematical Topics, 12,  Akademie Verlag, Berlin, 1997. 


\bibitem{Yag_Galst_CMP} 
K. Yagdjian and A. Galstian, {\it Fundamental solutions for the Klein-Gordon equation in de Sitter spacetime},
Comm. Math. Phys., 285 (2009), 293--344.

\bibitem{Macao} 
K. Yagdjian, {\it Integral transform approach to time--dependent partial differential equations},  Mathematical analysis, probability and applications--plenary lectures,  Springer Proc. Math. Stat., Vol. 177, Springer, (2016), 281--336.

\bibitem{JMP2019}
K. Yagdjian, {\it Global existence of the self--interacting scalar field in the de Sitter universe},  J. Math. Phys.,   60 (2019),  051503.

 
\bibitem{JDE2021}
K. Yagdjian and A. Galstian, {\it Fundamental solutions for the Dirac equation in curved spacetime and generalized Euler--Poisson--Darboux equation},  J. Differential Equations,  300 (2021), 80--117.
 

 \bibitem{ArXiv2023}
K. Yagdjian, {\it  Waves in cosmological background with static Schwarzschild radius in the expanding universe}, J. Math. Phys., 65, (2024), 083503. 
 
\end{thebibliography}
\end{document}